\documentclass[11pt,english,a4paper]{smfart}
\usepackage[english]{babel}
\usepackage{xcolor}
\usepackage{pbsi}
\usepackage[T1]{fontenc}
\usepackage{stmaryrd}
\usepackage{mathabx}
\usepackage{dsfont}
\usepackage[mathscr]{euscript}
\usepackage{cancel}
\usepackage{array}
\usepackage{float}
\usepackage{placeins}
\usepackage{csquotes}

\usepackage[pagebackref,colorlinks=true,linkcolor=blue,citecolor=blue,urlcolor=red, hypertexnames=true]{hyperref}
\usepackage[abbrev,backrefs,nobysame]{amsrefs}

\usepackage{caption}
 \usepackage{enumerate}

\usepackage{amsmath,amssymb,bm,amsfonts}
\usepackage{nicematrix}
\usepackage[all]{xy}
\entrymodifiers={+!!<0pt,\fontdimen22\textfont2>}

\usepackage{mathrsfs}

\usepackage{graphicx}

\usepackage{color}
\usepackage{comment}

\newtheorem{theorem}{Theorem}[section]
\newtheorem{lemma}[theorem]{Lemma}

\newtheorem{corollary}[theorem]{Corollary}
\newtheorem{proposition}[theorem]{Proposition}

\newcounter{hypothesis}
\renewcommand{\thehypothesis}{(H$_{\arabic{hypothesis}}$)}

\theoremstyle{remark}

\newtheorem*{remark*}{\it Remark}

\newcommand{\flba}[2]{
\xymatrix@C15pt{#1\ar@{|->}[r]&#2}}
\newcommand{\flcourte}[2]{
\xymatrix@C12pt{#1\ar[r]&#2}}

{\theoremstyle{definition}}

{\theoremstyle{definition}}

\makeatletter
\providecommand{\p@example}{\p@theorem}
\makeatother

\theoremstyle{definition}

{\theoremstyle{definition}}

\def\D{\ensuremath{\mathbb D}}

\def\T{\ensuremath{\mathbb T}}
\def\R{\ensuremath{\mathbb R}}
\def\Z{\ensuremath{\mathbb Z}}

\def\C{\ensuremath{\mathbb C}}

\def\N{\ensuremath{\mathbb N}}
\def\P{\ensuremath{\mathbb P}}
\def\E{\mathbb E}

\def\bth{\begin{theorem}}
\def\blm{\begin{lemma}}
\def\bpr{\begin{proposition}}
\def\bpf{\begin{proof}}
\def\epf{\end{proof}}
\def\epr{\end{proposition}}
\def\elm{\end{lemma}}
\def\eth{\end{theorem}}
\def\bco{\begin{corollary}}
\def\eco{\end{corollary}}
\def\be{\begin{enumerate}}
\def\ee{\end{enumerate}}
\def\bea{\begin{enumerate}[\rm (a)]}
\def\beun{\begin{enumerate}[\rm (1)]}
\def\bei{\begin{enumerate}[\rm (i)]}

\renewcommand{\Im}{\operatorname{Im}}
\renewcommand{\Re}{\operatorname{Re}}

\newcommand{\vertiii}[1]{{\left\vert\kern-0.25ex\left\vert\kern-0.25ex\left\vert #1 
    \right\vert\kern-0.25ex\right\vert\kern-0.25ex\right\vert}}

\newcommand{\lan}{\langle}
\newcommand{\ran}{\rangle}

\newcommand{\omg}{\omega}
\newcommand{\Omg}{\Omega}

\newcommand{\indic}{\mathds{1}}

\numberwithin{equation}{section}

\author[Valentin Gillet]{Valentin Gillet}
\address[Valentin Gillet]{Univ. Lille, CNRS, UMR 8524 - Laboratoire Paul Painlevé, F-59000 Lille, France}
\email{valentin.gillet@univ-lille.fr}

\thanks{This work was supported in part by the project COMOP of the French National Research Agency (grant ANR-24-CE40-0892-01). The author also acknowledges the support of the CDP C2EMPI, as well as the French State under the France-2030 programme, the University of Lille, the Initiative of Excellence of the University of Lille, and the European Metropolis of Lille for their funding and support of the R-CDP-24-004-C2EMPI project.}

\begin{document}

\title[Infinitely divisible invariant measures for linear operators on Banach spaces]{Infinitely divisible invariant measures for linear operators on Banach spaces}


\begin{abstract} 
We study the existence of stable and infinitely divisible invariant measures for operators on Banach spaces admitting a finite or countably infinite family of perfectly spanning unimodular eigenvector fields. These invariant measures are realized as the distributions of stochastic integrals with respect to stable and infinitely divisible complex random measures. We also study the existence of such invariant measures for semigroups of operators whose generators admit a finite or countably infinite family of perfectly spanning purely imaginary eigenvector fields.
\end{abstract}

\maketitle

\par\bigskip

\section{Introduction} \label{sectionIntroduction}

The aim of this paper is to complement the study \cite{Gill1} on the existence of infinitely divisible invariant measures for operators and semigroups of operators on infinite-dimensional separable complex Banach spaces. More precisely, Theorem 4.2 in \cite{Gill1} shows that, under certain geometric assumptions on the infinite-dimensional separable complex Banach space, any operator on such a Banach space admitting a sufficiently integrable family $(u_n)_{n \in J}$ of unimodular eigenvectors (with $J \subseteq \Z$) such that $\textrm{span} \{ u_n : n \in J \}$ is dense admits an isotropic\footnote{Isotropic means rotationally invariant} stable invariant measure with full support. The same result holds for semigroups of operators on such a Banach space for which there exists a sufficiently integrable family $(u_t)_{t \in J}$ (with $J \subseteq \R$) of purely imaginary eigenvectors for the generator of this semigroup such that $\textrm{span} \{ u_t : t \in J \setminus B \}$ is dense for every Borel subset $B$ of $J $ with Lebesgue measure zero (\cite[Thm. 3.8]{Gill1}). In the case of semigroups, if the Banach space is of type $2$, we obtain a strongly mixing isotropic Gaussian invariant measure with full support, whereas if the Banach space is of stable type $\alpha \in (0,2)$, we obtain an isotropic non-ergodic $\alpha$-stable invariant measure with full support. For operators, the $\alpha$-stable invariant measure constructed is never ergodic, even in the case where $\alpha = 2$ (i.e., the Gaussian case) by the result of Grivaux \cite[Thm. 5.1]{Grivaux11}. In order to obtain ergodic stable invariant measures with full support, we will need that the eigenvalues are somehow continuously distributed, exactly as in the Gaussian setting developed in the literature. 

\smallskip
The first to investigate which Hilbert spaces operators  are ergodic in the Gaussian sense was Flytzanis \cite{Flytzanis95}. Later on, Bayart and Grivaux continued this study in Banach spaces and also introduced the notion of a frequently hypercyclic operator (\cite{BayartGrivaux06} and \cite{BayartGrivaux07}). Let us recall their strategy to construct ergodic isotropic Gaussian invariant measures for operators.  In what follows, for every measure space $(J,\Sigma, m)$ and every $0 < p < \infty$, we denote by $L^p(J;E)$ the space of measurable functions $f : J \to E$ such that $\int_J \lVert f(t) \rVert^p m(dt) < \infty$. We usually equip $\T$ with the normalized Lebesgue measure but if $\nu$ is a Borel probability measure on $J = \T$, we rather write $L^p(\T;E;\nu)$ for this space instead of $L^p(\T;E)$. 

\bigskip
Let $E$ be an infinite-dimensional separable complex Banach space.
Let us recall that a finite or countably infinite sequence $(F_j)_{j \in I}$ of maps $F_j : \T \to E$ in $L^1(\T;E; \nu)$ such that $F_j(\lambda) \in \textrm{Ker}(T-\lambda)$ for each $\lambda \in \T$ is said to be $\nu$-spanning for an operator $T$ on $E$ if $\textrm{span}\{ F_j(\lambda) : j \in J, \lambda \in A \}$ is dense in $E$ for every Borel subset $A \subset \T$ such that $\nu(A) = 1$, where $\nu$ is a probability measure on $\T$. 
In the Gaussian setting, a main tool to construct ergodic isotropic Gaussian invariant measures with full support for operators admitting a perfectly spanning set of $\T$-eigenvectors is the covariance operator. Precisely, with each map $F : \T \to E$ in $L^2(\T;E;\nu)$, one can associate the operator $K_F : L^2(\T;\C;\nu) \to E$ defined by 
$$
K_F(f) = \int_{\T} f(\lambda) F(\lambda) d\nu(\lambda).
$$
When the Banach space $E$ is of type $p \in [1,2]$ and when the field $F$ is $\beta$-Hölderian for some $\beta > 1/p - 1/2$, then the operator $K_F$ is $\gamma$-radonifying. It implies that, when the operator $T$ on $E$ admits a finite or countably infinite sequence $(F_j)_{j \in I}$ of $\T$-eigenvector fields which is spanning with respect to the normalized Lebesgue measure $m$ on $\T$ and such that each $F_j$ is $\beta_j$-Hölderian for some $\beta_j > 1/p - 1/2$, then there exists a strongly mixing $T$-invariant isotropic Gaussian measure with full support on $E$. Such a Gaussian invariant measure is obtained via the construction of a covariance operator $R : E^* \to E$ using the intertwining equation \enquote{$TK = KM$}, where $K : \bigoplus\limits_{j\in I} L^2(\T;\C;m) \to E $ is defined for some family $(\alpha_j)_{j \in I}$ of positive numbers by
$$
K \left( \bigoplus\limits_{j\in I}f_j \right) = \sum_{j \in I} \alpha_j K_{F_j}(f_j),
$$
where $M : \bigoplus\limits_{j\in I} L^2(\T;\C;m) \to  \bigoplus\limits_{j\in I} L^2(\T;\C; m) $ 
is defined by  
$$
M(\bigoplus_{j \in I} f_j)(\lambda ) = \bigoplus_{j \in I} \lambda f_j(\lambda),
$$
and where the covariance operator is given by $R = K K^*$. We refer to \cite{BayartGrivaux06}, \cite{BayartGrivaux07}, \cite[Chapter 5]{BM09} and to \cite{BayartMatheron16} for more details.

\medskip
In the case of $\alpha$-stable measures, there is no covariance operator when the index of stability $\alpha$ is in $(0,2)$. A way to construct isotropic stable invariant measures for operators admitting perfectly spanning $\T$-eigenvector fields is to make use of the stochastic integral with respect to complex isotropic stable random measures, as constructed in \cite{Gill1}.
More precisely, it is shown in \cite[Thm. 4.1]{Gill1} that, if an operator on a complex separable Banach space of stable type $\alpha \in (0,2]$ has a dense and $\alpha$-integrable bilateral backward orbit, then it admits a strongly mixing isotropic $\alpha$-stable invariant measure with full support. A natural way to construct a bilateral backward orbit for an operator having an integrable $\T$-eigenvector field is to consider the sequence of the Fourier coefficients of this field. In the case where the operator has a finite or countably infinite sequence $(F_j)_{j \in I}$ of integrable $\T$-eigenvector fields, we will rather consider the stochastic integral $\int_{I \times \Z} \hat{F_j}(n) \, dM(j,n)$ with respect to a complex isotropic $\alpha$-stable random measure on $I \times \Z$ in order to take into account all the fields, where $\hat{F_j}(n)$ is the $n$-th Fourier coefficient of $F_j$.
Let us mention that the idea to consider random integrals to build invariant measures for operators admitting a dense bilateral backward orbit in \cite{Gill1} was inspired from the article \cite{AGNEESSENS_2026}, where random integrals with respect to Brownian motions are considered, but also from \cite{AGNEESSENS_2023}. The construction of a random integral with respect to random measures was inspired from \cite{Rosinski87}, \cite[Section 6.10]{Woyczynski19} (adapted to the complex setting in \cite{Gill1}), and to \cite[Chapters 3 and 6]{SamorodnitskyTaqqu94}.

\medskip
We will also study the existence of infinitely divisible ergodic invariant measure for strongly continuous semigroups of operators on Banach spaces to complement \cite[Thm. 3.8]{Gill1} using random integrals with respect to infinitely divisible complex random measures. In this direction, Chakir and El Mourchid \cite{ChakirMourchid18} seem to be the first to use stochastic integrals with respect to complex Brownian motions in order to construct strongly mixing Gaussian invariant measures with full support for semigroups of operators whose generators admit a large family of purely imaginary eigenvectors. Agneessens \cite{AGNEESSENS_2026} also obtained similar results to those of Chakir and El Mourchid. Let us also mention that, in order to construct ergodic Gaussian invariant measures for strongly continuous semigroups of operators for which there exists a finite or countably infinite family of purely imaginary perfectly spanning eigenvector fields for the generators, it is also possible to proceed exactly as in the case of operators, using the covariance operator (see \cite[Thm. 5.42]{BM09}). 

\smallskip
Finally, it is worth to note that any operator or semigroup of operators satisfying the so called Frequently Hypercyclic Criterion has a strongly mixing invariant measure with full support, see \cite{MurilloPeris13} for the operator case and \cite{MurilloPeris15} for the operator semigroup case.

\bigskip

We consider an infinite-dimensional separable complex Banach space $E$. If $\mu$ is a probability measure on $E$ equipped with its Borel subsets, we say that $\mu$ is infinitely divisible if for every integer $n \geq 1$, there exists a Borel probability measure $\mu_n$ on $E$ such that $\mu_n^{\star n} = \mu$.
Mau and Privault \cite{MauPrivault} studied conditions for an infinitely divisible invariant measure for an operator to be weakly mixing or strongly mixing with respect to this operator. Their paper also constitutes a motivation for this paper. An important class of infinitely divisible probability measures is given by stable measures.
Let us recall that an $E$-valued random variable $X$ is isotropic $\alpha$-stable with $\alpha \in (0,2]$ if for any $a,b >0$, the random variable $a X^{(1)} + b X^{(2)}$ has the same distribution as $(a^\alpha + b^\alpha)^{1/\alpha} X$, where $X^{(1)}$ and $X^{(2)}$ are independent copies of $X$, and if $\lambda X$ has the same distribution as $X$ for every complex number $\lambda$ such that $\lvert \lambda \rvert = 1$. Isotropic $2$-stable random vectors are Gaussian. It is well known that Gaussian probability measures have finite moments of all orders. When $\alpha <2$, an $\alpha$-stable probability measure on $E$ does not have a finite moment of order $\alpha$, but it has finite moments of order $p$ for every $0<p<\alpha$. All random variables and random integrals in this paper are defined on a suitable probability space $(\Omega,\mathcal F,\mathbb P)$.

\medskip

\bigskip
Let $(J, \Sigma, m)$ be a $\sigma$-finite measure space. Let us recall that a complex random measure on $(J, \Sigma, m)$ is a map $M$ on $\Sigma_0 := \{ A \in \Sigma : m(A) < \infty \}$ taking values in the space of complex random variables on $(\Omega,  \mathcal{F}, \P)$ such that for every pairwise disjoint sequence $(A_k)_{1 \leq k \leq n} \subset \Sigma_0$, the random vectors $(M(A_k))_{1 \leq k \leq n}$ are independent, and such that for every pairwise disjoint sequence $(A_k)_{k \geq 1} \subset \Sigma_0$ for which $\displaystyle \bigcup_{k \geq 1} A_k \in \Sigma_0$, we have that $\displaystyle \sum_{k \geq 1} M(A_k) = M \left(\displaystyle \bigcup_{k \geq 1} A_k \right) $ almost surely, where the series converges in probability (and hence almost surely). The independence of the complex-valued random variables $M(A_1), \dotsc , M(A_n)$ has to be understood in the sense that the vectors $ 
\begin{pmatrix}
M^{(1)}(A_1)\\
M^{(2)}(A_1)
\end{pmatrix}, \dotsc, 
\begin{pmatrix}
M^{(1)}(A_n)\\
M^{(2)}(A_n)
\end{pmatrix}$
are independent, where $M^{(1)}(A) := \Re(M(A))$ and $M^{(2)}(A) := \Im(M(A)) $. 
A standard way to construct a complex random measure with prescribed distribution is to apply Kolmogorov's existence theorem to a consistent family of finite-dimensional distributions. 

\medskip
We define random integrals with respect to random measures in a strong way. Thus, we need to assume certain geometric assumptions on the complex Banach space $E$. Let us recall that $E$ is of type $p \in [1,2]$ if for some (or, equivalently, every) finite $0<q < \infty$, there exists a constant $C_{p,q}$ such that, for all finite sequence $(x_n) \subset E$,
$$
\left(\E   \left \lVert \sum_n \varepsilon_n x_n \right \rVert^q  \right)^{1/q} \leq C_{p,q} \left(\sum_n \lVert x_n \rVert^p \right)^{1/p},
$$
where $(\varepsilon_n)_{n \geq 1}$ is a sequence of i.i.d Rademacher random variables. And $E$ is of stable type $p$ (with $0<p\leq 2$) if for some (or, equivalently, every) $0<q<p$ there exists a constant $C_{p,q}$ such that, for all finite sequence $(x_n) \subset E$,
$$
\left (\E \left \lVert \sum_n \theta_n^{(p)} x_n \right \rVert^q \right)^{1/q} \leq C_{p,q} \left(\sum_n \lVert x_n \rVert^p \right)^{1/p},
$$
where $( \theta_n^{(p)})_{n \geq 1}$ is a sequence of i.i.d standard symmetric $p$-stable random variables (i.e., real random variables satisfying $\E(e^{it \theta_n^{(p)}}) = \exp(-\lvert t \rvert^p)$ for every $t \in \R$). Let us recall also that if $E$ is of type $p \in [1,2]$ (resp. of stable type $p \in (0,2]$), then it is also of type $q \in [1,p]$ (resp. of stable type $q \in (0,p]$).
If $E$ is of stable type $p \in [1,2]$, then it is of type $p$, and if $E$ is of type $p \in [1,2]$, then $E$ is of stable type $q$ for any $0<q<p$, but not of stable type $p$ since $\theta_n^{(p)}$ does not have $p$-moment. Finally, $E$ is of type $2$ if and only if it is of stable type $2$. For more details on these notions, see \cite[p. 4]{Gill1}, \cite[Chapter 5]{LiQueffelec18}, \cite[Section 3.5]{Linde86}, \cite[Chapter III]{Schwartz81} and \cite[Sections 6.2 and 6.5]{Woyczynski19}. 

\medskip
Let us recall now the construction of the random integrals with respect to complex stable or infinitely divisible random measures. A complex random measure on a $\sigma$-finite measure space $(J,\Sigma, m)$ is a complex isotropic $\alpha$-stable random measure if for every $A \in \Sigma_0$, the real random variables $M^{(1)}(A) = \Re(M(A))$ and $M^{(2)}(A) = \Im(M(A))$ are jointly symmetric $\alpha$-stable with spectral measure $m(A) \times \gamma$, where $\gamma$ is the uniform probability measure on the unit sphere $\mathbb{S}_2$ of $\C$ (\cite[Section 6.1]{SamorodnitskyTaqqu94}). The measure $m$ is referred as the control measure of the random measure $M$. Suppose now that $E$ is a complex separable infinite-dimensional Banach space of stable type $\alpha \in (0;2]$. Let $M$ be a complex isotropic $\alpha$-stable random measure on the space $(J,\Sigma,m)$.
If $f = \displaystyle \sum_{j = 1}^n x_j \indic_{A_j}$ is a simple function, where the $A_j$'s are disjoint sets in $\Sigma_0$, then we set 
\begin{align*}
    \int_J f(t) M(dt) := \sum_{k=1}^n x_k M(A_k).
\end{align*}
For such a function $f$, the random integral $\int_J f(t) M(dt)$ is an isotropic $E$-valued $\alpha$-stable random vector. Moreover, for every $0 < q < \alpha$, there exists a constant $C>0$ depending on $q$ and $\alpha $ (but not on the function $f$) such that 
\begin{align} \label{inequalityqmomentstable}
    \left \lVert \int_J f(t) M(dt)  \right \rVert_q \leq C \left ( \int_J \lVert f(t) \rVert^\alpha m(dt) \right)^{1/\alpha},
\end{align}
see \cite[pages 7 and 10]{Gill1}.
Since the simple functions are dense in $L^\alpha(J;E)$, we can extend the definition of $\int_J f \, dM$ to every function $f \in L^\alpha(J;E)$. The random vector $\int_J f \, dM$ is thus isotropic $\alpha$-stable in $E$, satisfies (\ref{inequalityqmomentstable}), and its characteristic function is given by
\begin{align} \label{fonctioncaractvecteurstable}
    \E(e^{i \Re(\lan x^*, \int_J f \, dM \ran) }) = \exp \left( - c_0 \int_J \lvert \lan x^*, f(t) \ran \rvert^\alpha m(dt) \right),
\end{align}
where $c_0 := \frac{1}{2\pi} \int_0^{2\pi} \lvert \cos(x) \rvert^\alpha dx$. 

\smallskip
More generally, we can construct random integrals with respect to complex infinitely divisible random measures. Let $\mu$ be an infinitely divisible distribution $\C$, with characteristic function written for every $z \in \C$ as 
\begin{align*}
    \hat{\mu}(z) = \exp(\psi(z))
\end{align*}
with
\begin{align*}
    \psi(z) := -\frac{1}{2} q(z) + \int_{\C} (e^{i \Re(\overline{z}w)} -1 -i \Re(\overline{z} w) k(w)) \nu(dw),
\end{align*}
where $q : \C \to [0,\infty)$ is a quadratic form, $k$ is a truncation function and $\nu$ is a measure on $\C$ satisfying $\nu(\{ 0 \}) = 0$ and $\int_{\C} \min(\lvert w \rvert^2, 1) \nu(dw) < \infty$, see \cite[pp. 7--8]{Rosinski87}. A complex random measure $M$ on a $\sigma$-finite measure space $(J, \Sigma, m)$ is an infinitely divisible random measure with characteristics $[0,q,\nu]$ (written as $M \sim [0, q, \nu]$) if for every $A \in \Sigma$ with $m(A) < \infty$, we have that
\begin{align}
    \E(e^{i \Re(\overline{z} M(A))}) = \exp(m(A) \psi(z)), \; z \in \C. 
\end{align}
We define a random integral with respect to such a complex random measure $M$ in a strong way as for stable random measures. Let $E$ be a complex separable Banach space of type $p \in [1,2]$. We suppose that the infinitely divisible complex random measure $M$ also satisfies $\E(M(A)) = 0$ and that $\E(\lvert M(A) \rvert^p) \le C' m(A)$ for every $A \in \Sigma$ such that $m(A) < \infty$, where $C'>0$ is a constant independent of $A$. Then, as shown in \cite[pp. 24--25]{Gill1}, for every function $f \in L^p(J;E)$, the $E$-valued random vector $\int_J f \, dM$ has an infinitely divisible distribution with characteristic function
\begin{align}
    \E(e^{i \Re(\lan x^*, \int_J f \, dM \ran)}) = \exp \left (\int_J \psi(\overline{\lan x^*, f(t) \ran}) m(dt) \right),
\end{align}
and satisfies 
\begin{align} \label{inegalitéesperanceBanachtypep}
    \E \left \lVert \int_J f \, dM \right \rVert^p \le C \int_J \lVert f(t) \rVert^p m(dt)
\end{align}
for some constant $C>0$ depending on $p$ but not on $f$. Two important examples of such complex infinitely divisible random measures are the isotropic complex Gaussian random measures and the complex Poissonian random measures. 
Precisely, let $\nu$ be a symmetric complex probability measure on $\C$ such that $\nu(\{0 \}) = 0$ and $\int_{\C} \lvert w \rvert^2 \nu(dw) < \infty$, and let $\mu$ be the complex distribution satisfying $\hat{\mu}(z) = e^{\hat{\nu}(z) -1}$ for every $z \in \C$. Let $M$ be the complex random measure satisfying 
$$\E(e^{i \Re(\overline{z} M(A))}) = e^{m(A) (\hat{\nu}(z) -1)}$$ 
for every $z \in \C$. Let us notice that $\mu = \displaystyle \sum_{k \ge 0} \frac{\nu^k}{k!}$. The measure $\mu$ is called the exponent of the measure $\nu$ (see Linde \cite[Section 5.3]{Linde86}). Hence, for every measurable set $A$ with finite measure, the random variable $M(A)$ is the exponent of $m(A) \nu$. Moreover, we have that $\E(M(A)) = 0$, that $\E(\lvert M(A) \rvert^2) = m(A) \int_{\C} \lvert w \rvert^2 \nu(dw)$, and that $\textrm{supp}(\mu) = \overline{\displaystyle \bigcup_{k \ge 0} k \cdot \textrm{supp}(\nu)}$ (by \cite[Prop. 5.8.4]{Linde86}), showing that the random measure $M$ satisfies the necessary conditions to define the random integral when the Banach space $E$ is of type $2$. Suppose now that $M$ is a complex isotropic Gaussian random measure, with
$$
\E(e^{i \Re(\overline{z} M(A))}) = e^{-\frac{1}{2} \lvert z \rvert^2 m(A)}
$$
for every $z \in \C$. Then, we have that $\E(M(A)) = 0$ and that $\E(\lvert M(A) \rvert^2) = 2 m(A)$ for every measurable set $A$ with finite measure, showing again that we can define the random integral with respect to $M$ when $E$ is of type $2$. Of course, in the case where $M$ is a complex isotropic Gaussian random measure and when the Banach space $E$ is of type $2$, the construction was already done before since $E$ is of stable type $2$. We will see just after how to extend it to Banach spaces with type $p \in [1,2]$. 

\bigskip
We finish with two remarks. Let us suppose that $J = I \times \Z$, where $I$ is finite or countably infinite, and let us equip $J$ with its subsets and the counting measure. Then, first, for $f \in L^\alpha(J;E)$ or $f \in L^p(J;E)$ with $\alpha \in (0,2]$ and $p \in [1,2]$, we have that $$\int_J f \, dM = \displaystyle \sum_{(j,n) \in I \times \Z} M(\{(j,n) \}) f(j,n),$$ 
where $M$ is a complex isotropic $\alpha$-stable random measure on $J$ in the case where $E$ is of stable type $\alpha$ and $f \in L^\alpha(J;E)$, and $M$ is a complex infinitely divisible random measure on $J$ as above in the case where $E$ is of type $p$ and $f \in L^p(J;E)$. The convergence of the series holds almost surely in the following case: for every sequence $(J_k)_{k \ge 1}$ of nonempty finite subsets of $J$ such that $J_k \subseteq J_{k+1}$ for every $k \ge 1$ and $\displaystyle \bigcup_{k \ge 1} J_k = J$,
\begin{align*}
    \sum_{(j,n) \in I \times \Z} M(\{(j,n) \} f(j,n) = \lim_{k \to \infty} \sum_{(j,n) \in J_k} M(\{(j,n) \} f(j,n).
\end{align*}
Indeed, the definition of the random integral does not depend on the choice of the approximating sequence of step functions, and setting $f_k := f \indic_{J_k}$ for each $k \geq 1$, we have that $(f_k)_{k \ge 1}$ converges to $f$ in $L^\alpha(J;E)$ or in $L^p(J;E)$, depending on the case. Secondly, for every random variable $Y$ in $L^2(\Omg;E)$, we have that $\E(\lvert Y \rvert^p) \le 1 + \E(\lvert Y \rvert^2)$ for every $p \in [1,2]$. Suppose that $E$ is of type $p$ and suppose that $M$ is either a complex isotropic Gaussian random measure on $J$ or a complex Poissonian random measure on $J$. Then, since the measure of each set $\{(j,n) \}$ is equal to $1$ for every $(j,n) \in J$, we obtain that $\E(\lvert M(\{(j,n) \}) \rvert^p) \le C' m(\{(j,n) \})$ for every $(j,n) \in J$, for some constant $C'>0$. Thus, we can define the random integral with respect to $M$ even if $E$ is of type $p \in [1,2]$ and not only if $E $ is of type $2$, which is an improvement in the case of complex isotropic Gaussian random measures.

\section{Ergodic stable and infinitely divisible invariant measures for operators}

We study the existence of ergodic stable and infinitely divisible invariant measures for operators admitting a finite or countably infinite sequence of perfectly spanning $\T$-eigenvector fields. By \cite[Thm. 4.2]{Gill1} in the case where $0<\alpha <2$ and by \cite[Thm. 5.1]{Grivaux11} in the case where $\alpha = 2$, if an operator $T$ on a complex separable Banach space of stable type $\alpha$ admits a countable sequence $(u_n)_{n \in J}$ of unimodular $\alpha$-summable eigenvectors such that $\textrm{span}\{ u_n : n \in J \}$ is dense, then for every complex isotropic $\alpha$-stable random measure, the distribution of the vector $\displaystyle \sum_{n \in J} u_n M(\{n \})$ is a non-ergodic $\alpha$-stable $T$-invariant measure with full support. The aim of this section is to improve this result in the stable case as well as in the infinitely divisible case for operators admitting a perfectly spanning sequence of $\T$-eigenvector fields, as in the Gaussian case. The strategy is to construct a dense bilateral backward orbit obtained from each $\T$-eigenvector field of the operator. A natural way to do this is to consider the Fourier coefficients of each field. More precisely, if $T$ is an operator on $E$ and if $F : \T \to E$ is a $\T$-eigenvector field of $T$ in $L^1(\T;E)$, then denoting by $\hat{F}(n) = \int_{\T} F(z) \overline{z}^n dm(z)$ the $n$-th Fourier coefficient of $F$, we have that $T \hat{F}(n) = \int_{\T} z F(z) \overline{z}^n dm(z) = \hat{F}(n-1)$ for every $n \in \Z$. 

\subsection{General results}
The main results of existence of ergodic isotropic stable and infinitely divisible invariant measures for operators on an infinite-dimensional separable complex Banach space $E$ admitting perfectly spanning $\T$-eigenvector fields are the following. We start with the stable case.
\smallskip

\begin{theorem} \label{theoremstablechampsvpsunimod}
    Let $E$ be an infinite-dimensional separable complex Banach space of stable type $\alpha \in (0,2]$. Let $T$ be a bounded linear operator on $E$. Suppose that there exists a finite or countably infinite sequence $(F_j)_{j \in I}$ of maps in $L^1(\T;E)$ such that $T F_j(\lambda) = \lambda F_j(\lambda)$ for every $j \in I$ and $\lambda \in \T$, and that for every Borel subset $A \subseteq \T$ such that $m(A) = 1$, the subspace $\textrm{span}\{ F_j (\lambda) : j \in I, \lambda \in A \}$ is dense in $E$. Suppose that for every $j \in I $, the series $\displaystyle \sum_{n \in \Z} \lVert \hat{F_j}(n) \rVert^\alpha$ is convergent. Then $T$ admits a strongly mixing isotropic $\alpha$-stable invariant measure with full support.
\end{theorem}

\smallskip

\begin{proof}
Let us consider a complex isotropic $\alpha$-stable random measure $M$ on $I \times \Z$ equipped with its subsets and the counting measure, which is $\sigma$-finite. For each $j \in I$ and $n \in \Z$, let $x_{j,n} := \hat{F_j}(n) = \int_{\T} \overline{\lambda}^n F_j(\lambda) dm(\lambda)$. Then $T x_{j,n} = x_{j, n-1}$ for each $j \in I$ and $n \in \Z$. 
    \smallskip

    Let us remark that the subspace $\textrm{span}\{ x_{j,n} : j \in I, n \in \Z \}$ is dense in $E$. Indeed, if $x^* \in E^*$ vanishes on $\textrm{span}\{ x_{j,n} : j \in I, n \in \Z \}$, then $\int_{\T} \overline{\lambda}^n \lan x^*, F_j(\lambda) \ran \, dm(\lambda) = 0 $ for every $j \in I$ and $n \in \Z$. Since the functions $\lambda \mapsto \lan x^*, F_j(\lambda) \ran$ belong to $L^1(\T;\C)$, this implies that for every $j \in I$, there exists a Borel subset $A_j \subseteq \T$ with $m(A_j) = 1$ such that $\lan x^*, F_j(\lambda) \ran = 0$ for every $\lambda \in A_j$. Thus, since $I$ is finite or countably infinite, the Borel subset $A := \displaystyle \bigcap_{j \in I} A_j$ has full measure and for every $j \in I$ and $\lambda \in A$, we have that $\lan x^*, F_j(\lambda) \ran = 0$. This implies that $x^*$ vanishes on $\textrm{span}\{ F_j(\lambda) : j \in I, \lambda \in A \}$ and so $x^* = 0$.

    \smallskip
    Now let $(a_j)_{j \in I}$ be a sequence of positive real numbers such that $$\displaystyle \sum_{(j,n) \in I \times \Z} a_j^\alpha \lVert x_{j,n} \rVert^\alpha < \infty.$$
    Then the random vector 
    \begin{align*}
        G := \int_{I \times \Z} a_j \, x_{j,n} \, dM(j,n) = \sum_{(j,n) \in I \times \Z} a_j x_{j,n} M_{j,n},
    \end{align*}
    where $M_{j,n} := M(\{(j,n)\})$ for every $(j,n) \in I \times \Z $,
    is almost surely well-defined on $(\Omega, \mathcal{F}, \P)$ and is isotropic $\alpha$-stable. Moreover, using the characteristic function (\ref{fonctioncaractvecteurstable}), the distribution $\mu_G$ of the vector $G$ is easily seen to be $T$-invariant since
    \begin{align*}
        \hat{\mu_G}(T^* x^*) = \exp \left(-c_0 \sum_{(j,n) \in I \times \Z} a_j^\alpha\lvert  \lan x^*, x_{j,n-1} \ran \rvert^\alpha \right) = \hat{\mu_G}(x^*)
    \end{align*}
    for every $x^* \in E^*$. Since the subspace $\textrm{span}\{ x_{j,n} : j \in I, n \in \Z \}$ is dense in $E$ and since the support of $\mu_G$ contains $\{ x^* \in E^* : \hat{\mu_G}(t x^*) = 1 , \forall t \in \R \}^{\perp}$ by \cite[Lemma 3.6]{Gill1}, we see that the random vector $G$ has full support.

    \medskip
    It only remains to prove that the distribution of $G$ is strongly mixing with respect to $T$. It is possible to just notice that the proof of mixing in \cite[Thm 3.1]{Gill1} can be adapted to our context. Instead, we provide an independent proof. We will conjugate the dynamical system $(E, \textrm{Bor}(E), \mu_G, T)$ to a Bernoulli shift, since $T$ acts as a backward shift on the second coordinate. Let us set $\Omega_0 := \C^{I \times \Z}$, $\mathcal{F}_0 := \bigotimes\limits_{(j, n)\in I \times \Z} \mathcal{\textrm{Bor}}(\C)$ and $\P_0 := \bigotimes\limits_{(j, n)\in I \times \Z} \nu_{\alpha}$, where $\nu_\alpha$ denotes the common distribution of the complex random variables $M_{j,n}$. Equivalently, $\mathcal{F}_0$ is the $\sigma$-algebra generated by the cylinders of $\Omega_0$. Denote by $\psi : \Omg \to \Omg_0 $ the map defined by $\psi(\omg) = (M_{j,n}(\omg))_{(j,n) \in I \times \Z}$. Then we have that $\P \circ \psi^{-1} = \P_0$. Let $S_0 : \Omg_0 \to \Omg_0$ be the forward shift defined by $S_0( (\omg_{j,n})_{(j,n) \in I \times \Z}) = (\omg_{j,n+1})_{(j,n) \in I \times \Z} $, and let $\Phi : \omg = (\omg_{j,n} )_{(j,n) \in I \times \Z} \in \Omg_0 \mapsto \displaystyle \sum_{(j,n) \in I \times \Z} a_j \, x_{j,n} \, \omg_{j,n}$. Let us notice that the series defining $\Phi(\omega)$ converges $\P_0$-almost surely, since the series defining $G(\omega)$ converges $\P$-almost surely and since $\P \circ \psi^{-1} = \P_0$. Then we notice that $\Phi \circ S_0 = T \circ \Phi$, and since $G = \Phi \circ \psi$, we also notice that $\P_0 \circ \Phi^{-1} = \mu_G$. This means that $(E, \textrm{Bor}(E), \mu_G, T)$ is a factor of $(\Omg_0, \mathcal{F}_0, \P_0, S_0)$. Thus, in order to show that $\mu_G$ is strongly mixing with respect to $T$, it is enough to show that the dynamical system $(\Omg_0, \mathcal{F}_0, \P_0, S_0)$ is strongly mixing.
    \smallskip
    
    In order to see this, let us consider the dynamical system $(\Tilde{\Omg_0}, \Tilde{\mathcal{F}_0}, \Tilde{\P_0}, \Tilde{S_0})$, where $\Tilde{\Omg_0} := (\C^{I})^{\Z}$, $\Tilde{\mathcal{F}_0} = \bigotimes\limits_{n \in \Z} (\bigotimes\limits_{j\in I } \textrm{Bor}(\C^I))$, $\Tilde{\P_0} = \bigotimes\limits_{n \in \Z} (\bigotimes\limits_{j\in I } \nu_\alpha) $ and $\Tilde{S_0} : (\omg^{(n)})_{n \in \Z} \mapsto (\omg^{(n+1)})_{n \in \Z}$ for each $\omg^{(n)} = (\omg_{j,n})_{j \in I} \in \C^I$. Then $\Tilde{\P_0} = \P_0 \circ R^{-1}$ and $R \circ S_0 = \Tilde{S_0} \circ R$, where $R : \Omg_0 \to \Tilde{\Omg_0}$ is invertible and given by $R((\omg_{j,n})_{(j,n) \in I \times \Z}) = (\omg^{(n)})_{n \in \Z}$ with $\omg^{(n)} = (\omg_{j,n})_{j \in I}$ for each $n \in \Z$. Hence the dynamical system $(\Omg_0, \mathcal{F}_0, \P_0, S_0)$ is isomorphic to $(\Tilde{\Omg_0}, \Tilde{\mathcal{F}_0}, \Tilde{\P_0}, \Tilde{S_0})$, and is strongly mixing (see, for instance, \cite[Examples 4.2.3 and 7.1.3]{VianaOliveira16}).

    \smallskip
    This concludes the proof of Theorem \ref{theoremstablechampsvpsunimod}.
    \end{proof}

\smallskip
        
Theorem \ref{theoremstablechampsvpsunimod} is stated for the normalized Lebesgue measure on $\T$, which is the most common situation. We can also cite it for a $\nu$-perfectly spanning sequence of $\T$-eigenvector fields for $T$, when $\nu$ is a probability measure on $\T$.   

\smallskip

\bth
Let $E$ be an infinite-dimensional separable complex Banach space of stable type $\alpha \in (0,2]$. Let $T$ be a bounded linear operator on $E$. Suppose that there exists a finite or countably infinite sequence $(F_j)_{j \in I}$ of maps $F_j : \T \to E$ satisfying $\int_{\T} \lVert F_j(\lambda) \rVert \nu(d\lambda) < \infty$ for each $j \in I$ such that $T F_j(\lambda) = \lambda F_j(\lambda)$ for every $j \in I$ and $\lambda \in \T$, and that for every Borel subset $A \subseteq \T$ such that $\nu(A) = 1$, the subspace $\textrm{span}\{ F_j (\lambda) : j \in I, \lambda \in A \}$ is dense in $E$, where $\nu$ is a probability measure on $\T$. Suppose that for every $j \in I $, the series $\displaystyle \sum_{n \in \Z} \lVert x_{j,n} \rVert^\alpha$ is convergent, where
$$
x_{j,n} := \int_{\T} \overline{\lambda}^n F_j(\lambda) \, \nu(d\lambda)
$$
for every $n \in \Z$ and $j \in I$.
Then $T$ admits a strongly mixing isotropic $\alpha$-stable invariant measure with full support.
\eth

\smallskip
\bpf
Exactly as in the proof of Theorem \ref{theoremstablechampsvpsunimod}, we have that $T x_{j,n} = x_{j, n-1}$ for every $n \in \Z$ and $j \in I$. Moreover, the subspace $\textrm{span}\{ x_{j,n} : j \in I, n \in \Z \}$ is dense in $E$. Indeed, if a linear functional $x^* \in E^*$ vanishes on this subspace, then for every $n \in \Z$, we have that $\int_{\T} \overline{\lambda}^n \lan x^*, F_j(\lambda) \ran \nu(d\lambda) = 0$. This means that the Fourier coefficients of the measure $\lan x^*, F_j(\cdot) \ran \nu$ are all equal to $0$, and thus $\lan x^*, F_j(\cdot) \ran = 0$ in $L^1(\T;\C;\nu)$. The rest of the proof is the same to that one of Theorem \ref{theoremstablechampsvpsunimod}. 
\epf

\smallskip
The analogue of Theorem \ref{theoremstablechampsvpsunimod} for infinitely divisible invariant measures is the following.

\smallskip

\bth \label{theoremIDchampsvpsunimod}
Let $E$ be an infinite-dimensional separable complex Banach space of type $p \in [1,2]$. Let $T$ be a bounded linear operator on $E$. Suppose that there exists a finite or countably infinite sequence $(F_j)_{j \in I}$ of maps $F_j : \T \to E$ satisfying $\int_{\T} \lVert F_j(\lambda) \rVert \pi(d\lambda) < \infty$ for each $j \in I$ such that $T F_j(\lambda) = \lambda F_j(\lambda)$ for every $j \in I$ and $\lambda \in \T$, and that for every Borel subset $A \subseteq \T$ such that $\pi(A) = 1$, the subspace $\textrm{span}\{ F_j (\lambda) : j \in I, \lambda \in A \}$ is dense in $E$, where $\pi$ is a probability measure on $\T$. Suppose that for every $j \in I $, the series $\displaystyle \sum_{n \in \Z} \lVert x_{j,n} \rVert^p$ is convergent, where
$$
x_{j,n} := \int_{\T} \overline{\lambda}^n F_j(\lambda) \, \pi(d\lambda)
$$
for every $n \in \Z$ and $j \in I$.
Then for every complex infinitely divisible random measure $M \sim [0, q, \nu]$ on $I \times \Z$ satisfying $\E(M(A)) = 0$ and $\E(\lvert M(A) \rvert^p) \le C' \#(A)$ for every finite subset $A$ of $I \times \Z$ and for some constant $C'>0$, where $\# $ is the counting measure on $I \times \Z$, the infinitely-divisible distribution of 
\begin{align*}
    \sum_{(j,n) \in I \times \Z} a_j x_{j,n} M(\{(j,n) \})
\end{align*}
is $T$-invariant and strongly mixing, where $(a_j)_{j \in I}$ is a sequence of positive real numbers such that $\displaystyle \sum_{(j,n) \in I \times \Z} a_j^p \lVert x_{j,n} \rVert^p < \infty$. If moreover the distribution of the $M(\{(j,n)\})$'s has full support, then so does the infinitely divisible $T$-invariant measure constructed.
\eth

\smallskip
\bpf
We can suppose that $I \subseteq \N$.
The proof of Theorem \ref{theoremIDchampsvpsunimod} is essentially analogous to the proof of Theorem \ref{theoremstablechampsvpsunimod}. Precisely, in order to prove that the distribution $\mu_G$ of the random vector $G := \int_{I \times \Z} a_j x_{j,n} \, dM(j,n)$ is $T$-invariant, we make use of the characteristic function 
\begin{align*}
    \hat{\mu_G}(x^*) = \exp \left(\sum_{(j,n) \in I \times \Z} \psi(\overline{\lan x^*, a_j x_{j,n} \ran}) \right),
\end{align*}
together with the fact that $\displaystyle \sum_{(j,n) \in I \times \Z} \lvert  \psi(\overline{\lan x^*, a_j x_{j,n} \ran}) \rvert < \infty$ by Rosinski \cite[Lemma 3.1.4]{Rosinski87}.

It only remains to show that if the common distribution of the random variables $M(\{(j,n) \})$ has full support, then the distribution of the random vector 
$$ G = \sum_{(j,n) \in I \times \Z} a_j x_{j,n} M(\{(j,n) \})$$ 
has full support too. However, this follows from an easy adaptation of the proof given in \cite[Thm. 3.9]{AGNEESSENS_2023}, using the fact that $\E(\lvert M(\{(j,n) \}) \rvert^p) \leq C'$ for every $(j,n) \in I \times \Z$, as well as by \cite[Prop. 6.7.1]{Woyczynski19}. 
\epf

\smallskip
In particular, Theorem \ref{theoremIDchampsvpsunimod} gives the following version for isotropic Gaussian measures, improving a bit Theorem \ref{theoremstablechampsvpsunimod}.
\smallskip

\bco \label{corexistencemesureGaussinvtypep}
Let $E$ be an infinite-dimensional separable complex Banach space of type $p \in [1,2]$. Let $T$ be a bounded linear operator on $E$. Suppose that there exists a finite or countably infinite sequence $(F_j)_{j \in I}$ of maps in $L^1(\T;E)$ such that $T F_j(\lambda) = \lambda F_j(\lambda)$ for every $j \in I$ and $\lambda \in \T$, and that for every Borel subset $A \subseteq \T$ such that $m(A) = 1$, the subspace $\textrm{span}\{ F_j (\lambda) : j \in I, \lambda \in A \}$ is dense in $E$. Suppose that for every $j \in I $, the series $\displaystyle \sum_{n \in \Z} \lVert \hat{F_j}(n) \rVert^p$ is convergent. Then $T$ admits a strongly mixing isotropic Gaussian invariant measure with full support.
\eco

\smallskip
For Banach spaces with a finite cotype, when the operator $T$ admits a $\beta$-Hölderian perfectly spanning $\T$-eigenvector field for some $\beta>1/2$, we have the following result for the existence of a strongly mixing invariant measure. 

\smallskip

\bpr \label{propBanachcotypemesureinvariante}
Let $E$ be a complex infinite-dimensional separable Banach space with a finite cotype. Let $T$ be a bounded linear operator on $E$ and suppose that there exists a map $F : \T \to E $ such  that $T F(\lambda) = \lambda F(\lambda)$ for every $\lambda \in \T$. Suppose that the map $F$ is $\beta$-Hölderian for some $\beta>1/2$. Then for every sequence of centered i.i.d complex random variables $X_n$, $n \in \Z$, with finite moments of any order, the distribution of the random vector $G :=\displaystyle \sum_{n \in \Z} X_n \hat{F}(n)$
is $T$-invariant and strongly mixing. Suppose moreover that, for every Borel subset $A \subseteq \T$ such that $m(A) = 1$, the subspace $\textrm{span}\{ F (\lambda) : \lambda \in A \}$ is dense in $E$, and that the random variable $X_0$ has full support. Then the distribution of $G$ has full support. 
\epr
\smallskip
\bpf
The proof of the first part of Proposition \ref{propBanachcotypemesureinvariante} is based on \cite[Remark 3.5]{BayartGrivaux07}, as well as on a theorem of Maurey and Pisier (\cite[Thm. 5.4]{ImekrazRobertThomann16}).
Note that for each $x^* \in E^*$, the scalar map $\lan x^*, F(.) \ran$ is $\beta$-Hölderian with $\beta>1/2$. By a theorem of Bernstein, the series $\displaystyle \sum_{n \in \Z} \lvert \lan x^*, \hat{F}(n) \ran \rvert$ is thus convergent. Since $E$ has a finite cotype and thus does not contain a copy of $c_0$, the series $\displaystyle \sum_{n \in \Z} \hat{F}(n) $ converges unconditionally in $E$ by a theorem of Bessaga and Pe{\l}czy{\'n}ski (see \cite[p. 22]{DiestelJarchowTonge95}). In particular, if $(\varepsilon_n)_{n \in \Z}$ is a sequence of i.i.d Rademacher random variables, then  both the series $\displaystyle \sum_{n \ge 1} \varepsilon_n \hat{F}(n)$ and $\displaystyle \sum_{n \le -1} \varepsilon_n \hat{F}(n)$ are convergent almost surely. By \cite[Thm. 5.4]{ImekrazRobertThomann16}, this implies that both the series $\displaystyle \sum_{n \ge 1} X_n \hat{F}(n)$ and $\displaystyle \sum_{n \le  -1} X_n \hat{F}(n)$ are convergent almost surely too. This shows that the random vector $G$ is well-defined. Moreover, it is not difficult to show that the distribution of $G$ is $T$-invariant, and the proof of mixing given in Theorem \ref{theoremstablechampsvpsunimod} shows that this distribution is strongly mixing. It remains to prove that the distribution of $G$ has full support.

Let $y \in E$ and $\varepsilon >0$. We have to prove that $\P(\lVert G - y \rVert < \varepsilon) >0$. We proceed as in \cite[Proof of Thm. 1.2]{MouzeMunnier21}. Let us note that $\textrm{span}\{ \hat{F}(n) : n \in \Z \}$ is dense in $E$. Thus, we can suppose that $y = \displaystyle \sum_{n=-d}^{d} a_n \hat{F}(n)$ for some integer $d \ge 1$ and $a_n \in \C$. Since $\displaystyle \sum_{\lvert k \rvert \le N} X_k \hat{F}(k)$ converges almost surely to the random vector $G$, we have that $\displaystyle \sum_{\lvert k \rvert > N} X_k \hat{F}(k)$ converges to $0$ in probability. Hence, we can choose an integer $N > d$ such that $\P \left(\lVert \displaystyle \sum_{\lvert k \rvert > N} X_k \hat{F}(k) \rVert < \varepsilon/2 \right) >0$. Let us set $a_n = 0$ for $d < \lvert n \rvert \le N$. Then
\smallskip
\begin{align*}
    &\P(\lVert G - y \rVert < \varepsilon) \\ 
    &\geq \P \left(\lVert \sum_{\lvert k \rvert \le N} X_k \hat{F}(k) - y \rVert < \varepsilon/2 \quad \textrm{and} \quad \lVert \sum_{\lvert k \rvert > N} X_k \hat{F}(k) \rVert < \varepsilon/2 \right) \\
    &\ge \P \left(\lVert \sum_{\lvert k \rvert \le N} X_k \hat{F}(k) - y \rVert < \varepsilon/2 \right) \cdot \P\left(\lVert \sum_{\lvert k \rvert > N} X_k \hat{F}(k) \rVert < \varepsilon/2 \right) \\
    &\ge \prod_{\lvert k \rvert \le N} \P \left(\lvert X_k - a_k \rvert \lVert \hat{F}(k) \rVert < \frac{\varepsilon}{2N+1} \right) \cdot \P\left(\lVert \sum_{\lvert k \rvert > N} X_k \hat{F}(k) \rVert < \varepsilon/2 \right) \\
    &>0.
\end{align*}
Here, we used the independence of the random variables $X_n$, $n \in \Z$, as well as the fact that these random variables have full support. This concludes the proof of Proposition \ref{propBanachcotypemesureinvariante}.
\epf

\medskip

We now discuss sufficient conditions in order to satisfy Theorems \ref{theoremstablechampsvpsunimod} and \ref{theoremIDchampsvpsunimod}, first in the case of a Hilbert space, which is the easiest case, and next in the case of a Banach space.

\subsection{Sufficient conditions in a Hilbert space}

In this subsection, we denote by $H$ an infinite-dimensional separable complex Hilbert space. Note that $H$ is of stable type $2$ and type $2$. Moreover, when the circle $\T$ is equipped with the normalized Lebesgue measure $m$, then Parseval's equality holds in $L^2(\T;H)$: for every $f \in L^2(\T;H)$, we have that
\begin{align*}
    \sum_{n \in \Z} \lVert \hat{f}(n) \rVert^2 = \int_{\T} \lVert f(\lambda) \rVert^2 dm(\lambda).
\end{align*}
Indeed, let $(e_j)_{j \geq 1}$ be an orthonormal basis of $H$. Then $f(\lambda) = \displaystyle\sum_{j \geq 1} f_j(\lambda) e_j$, and so $\lVert f \rVert_{L^2(\T;H)}^2 = \displaystyle \sum_{j \geq 1} \lVert f_j \rVert_{L^2(\T;\C)}^2$. Thus, the result follows from the scalar case. We immediately deduce the following result regarding Gaussian invariant measures. This recovers the result of Bayart and Grivaux \cite[Thm. 5.38]{BM09} in a Hilbert space for the normalized Lebesgue measure.

\smallskip

\bpr
Let $H$ be an infinite-dimensional separable complex Hilbert space and let $T$ be a bounded linear operator on $H$. Suppose that there exists a finite or countably infinite sequence $(F_j)_{j \in I}$ of maps in $L^2(\T;H)$ such that $T F_j(\lambda) = \lambda F_j(\lambda)$ for every $j \in I$ and $\lambda \in \T$, and that for every Borel subset $A \subseteq \T$ such that $m(A) = 1$, the subspace $\textrm{span}\{ F_j (\lambda) : j \in I, \lambda \in A \}$ is dense in $H$. Then the operator $T$ admits a strongly mixing isotropic Gaussian invariant measure with full support.
\epr

\smallskip
We also have the following analogue for infinitely divisible distributions.

\smallskip
\bpr
Let $H$ be an infinite-dimensional separable complex Hilbert space and let $T$ be a bounded linear operator on $H$. Suppose that there exists a finite or countably infinite sequence $(F_j)_{j \in I}$ of maps in $L^2(\T;H)$ such that $T F_j(\lambda) = \lambda F_j(\lambda)$ for every $j \in I$ and $\lambda \in \T$, and that for every Borel subset $A \subseteq \T$ such that $m(A) = 1$, the subspace $\textrm{span}\{ F_j (\lambda) : j \in I, \lambda \in A \}$ is dense in $H$. Then for every complex infinitely divisible random measure $M \sim [0, q, \nu]$ on $I \times \Z$ satisfying $\E(M(A)) = 0$ and $\E(\lvert M(A) \rvert^2) \le C' \#(A)$ for every finite subset $A$ of $I \times \Z$ and for some constant $C'>0$, where $\#$ is the counting measure on $I \times \Z$, the infinitely-divisible distribution of 
\begin{align*}
    \sum_{(j,n) \in I \times \Z} a_j \hat{F_j}(n) M(\{(j,n) \})
\end{align*}
is $T$-invariant and strongly mixing , where $(a_j)_{j \in I}$ is a sequence of positive real numbers such that $\displaystyle \sum_{(j,n) \in I \times \Z} a_j^2 \lVert \hat{F_j}(n) \rVert^2 < \infty$.
\epr

Let us note that, in general, the Fourier coefficients of a function $f \in L^2(\T;H)$ may not be always easy to compute. For this reason, it is more convenient to work with the $\T$-eigenvector fields instead of the Fourier coefficients of these fields. We now introduce two function spaces and will show that these function spaces coincide for $s \in (0,1)$. This is inspired by \cite[Prop. 1.4]{BenyiOh13} and \cite[p. 240]{Hormander90}.
 
\smallskip 
 For every $s >0$, we denote by $\mathcal{H}^s(\T;H)$ the space of functions $f \in L^2(\T;H)$ such that $\displaystyle \sum_{n \in \Z} (1 + \lvert n \rvert^2)^s \lVert \hat{f}(n) \rVert^2 < \infty$. Then we have the following useful lemma.

\smallskip

\blm
Let $\alpha \in (0,2]$. Suppose that there exists $s > 1/\alpha - 1/2$ such that $f \in \mathcal{H}^s(\T;H)$. Then $\displaystyle \sum_{n \in \Z} \lVert \hat{f}(n) \rVert^\alpha < \infty $.
\elm

\smallskip

\bpf
This follows from Hölder's inequality. Indeed, let $p = 2/\alpha$ and let $q = 2/(2-\alpha)$ the conjugate exponent of $p$\footnote{If $\alpha = 2$, then we just set $p = 1$ and $q = \infty$.}. Then
\begin{align*}
    \sum_{n \ne 0} \lVert \hat{f}(n) \rVert^\alpha &\leq (\sum_{n \ne 0} \lvert n \rvert^{2s} \lVert \hat{f}(n) \rVert^2)^{\alpha/2} (\sum_{n \ne 0} \lvert n \rvert^{-s \alpha q})^{1/q} \\
    &\leq (\sum_{n \ne 0} \lvert n \rvert^{2s} \lVert \hat{f}(n) \rVert^2)^{\alpha/2} (\sum_{n \ne 0} \lvert n \rvert^{-\frac{2 s \alpha}{2-\alpha}})^{(2-\alpha)/2},
\end{align*}
which is finite since $f \in \mathcal{H}^s(\T;H)$ and since $2s \alpha/(2-\alpha) > 1$.
\epf

\smallskip
For $s \in (0,1)$, let $W^{s,2}(\T;H) $ be the space of functions $f \in L^2(\T;H)$ such that $$[f]_{W^{s,2}}^2 := \int \int_{(-\pi, \pi)^2} \frac{\lVert f(e^{it}) - f(e^{iu}) \rVert^2 }{\lvert e^{it} - e^{i u} \rvert^{2 s +1}} dt du < \infty.$$

\smallskip

\bpr \label{propégalitéespacesSobolev}
For every $s \in (0,1)$, we have that $\mathcal{H}^s(\T;H) = W^{s,2}(\T;H)$. 
\epr

\smallskip
\bpf
First, note that 
\begin{align*}
    \int \int_{(-\pi, \pi)^2} \frac{\lVert f(e^{it}) - f(e^{iu}) \rVert^2 }{\lvert e^{it} - e^{i u} \rvert^{2 s +1}} dt du &= 
    \int \int_{(-\pi, \pi)^2} \frac{\lVert f(e^{i(u+h)}) - f(e^{iu}) \rVert^2 }{\lvert e^{i(u+h)} - e^{i u} \rvert^{2 s +1}} dh du \\ &= \int \int_{(-\pi, \pi)^2} \frac{\lVert f(e^{i(u+h)}) - f(e^{iu}) \rVert^2 }{\lvert e^{ih} - 1 \rvert^{2 s +1}} dh du.
\end{align*}
Moreover, by Parseval's identity, we have that
\begin{align*}
    \int_{-\pi}^{\pi} \lVert f(e^{i(h+u)}) - f(e^{iu}) \rVert^2 du = 2 \pi \sum_{n \in \Z} \lvert e^{inh} -1 \rvert^2 \lVert \hat{f}(n) \rVert^2
\end{align*}
for every $h \in (-\pi,\pi)$. Thus $[f]_{W^{s,2}}^2 = 2\pi \displaystyle \sum_{n \in \Z} \lVert \hat{f}(n) \rVert^2 I_n$, where
\begin{align*}
    I_n := \int_{-\pi}^{\pi} \frac{\lvert e^{inh} -1 \rvert^2}{\lvert e^{ih} -1 \rvert^{2s +1}} dh.
\end{align*}
Note that $I_0 = 0$, $I_{-n} = I_n$ for every $n \geq 1$ and that $I_n = 2 \int_0^{\pi}  \frac{\lvert e^{inh} -1 \rvert^2}{\lvert e^{ih} -1 \rvert^{2s +1}} dh $ for every $n \geq 1$. Using that $\frac{2}{\pi} h \leq \lvert e^{ih} -1 \rvert \leq h$ for every $h \in (0,\pi)$, we thus have for $n \geq 1$,
\begin{align*}
    I_n &\lesssim  \int_0^{\pi} \frac{\min(4, n^2 h^2)}{h^{1+2s}} dh \\
    &\lesssim (\int_0^{1/n} n^2 h^{1-2s} dh +  \int_{1/n}^\pi h^{-1-2s} dh) \\
    &\lesssim \frac{1}{2-2s} n^{2s} + \frac{n^{2s}}{2s} \\
    &\lesssim \, n^{2s}.
\end{align*}
Moreover, for $h \in (0,1/n]$, we have $\frac{\lvert e^{inh} -1 \rvert^2}{\lvert e^{ih} -1 \rvert^{1 + 2s}} \geq \frac{4 n^2 h^2}{\pi^2 h^{1+2s}} \geq \frac{4}{\pi^2} n^2 h^{1-2s}$, and thus
\begin{align*}
    I_n \gtrsim \int_0^{1/n} n^2 h^{1-2s} dh \gtrsim \, n^{2s}.
\end{align*}
It follows that $I_n \asymp \lvert n \rvert^{2s}$ for all $n \in \Z$. 

\smallskip
Finally, we obtain that
\begin{align*}
    [f]_{W^{s,2}}^2 &= 2 \pi \sum_{n \in \Z} I_n \, \lVert \hat{f}(n) \rVert^2 \\
    & \asymp \sum_{n \in \Z} \lvert n \rvert^{2s} \lVert \hat{f}(n) \rVert^2 \\
    & \asymp \sum_{n \in \Z} (1 + \lvert n \rvert^2)^s \lVert \hat{f}(n) \rVert^2.
\end{align*}
This concludes the proof of Proposition \ref{propégalitéespacesSobolev}.
\epf

\smallskip
We immediately see that every $\beta$-Hölderian function $f : \T \to H$ belongs to $\mathcal{H}^s(\T;H)$ if $0 < s < \beta <1$. We thus have the following result for stable invariant distributions.

\smallskip

\bco \label{cormesuresstablesinvholdercashilbert}
Let $H$ be an infinite-dimensional separable complex Hilbert space and let $T$ be a bounded linear operator on $H$. Let $\alpha \in (2/3,2]$. Suppose that there exists a finite or countably infinite sequence $(F_j)_{j \in I}$ of maps $F_j : \T \to H$ such that $T F_j(\lambda) = \lambda F_j(\lambda)$ for every $j \in I$ and $\lambda \in \T$, and that for every Borel subset $A \subseteq \T$ such that $m(A) = 1$, the subspace $\textrm{span}\{ F_j (\lambda) : j \in I, \lambda \in A \}$ is dense in $H$. Suppose also that for each $j \in I$, the map $F_j$ is $\beta_j$-Hölderian for some $\beta_j > 1/\alpha- 1/2$. Then the operator $T$ admits a strongly mixing isotropic $\alpha$-stable invariant measure with full support.
\eco

\medskip
Similarly, we obtain the following result for infinitely divisible distributions, since $H$ is of type $p \in [1,2]$.

\medskip
\bco
Let $H$ be an infinite-dimensional separable complex Hilbert space and let $T$ be a bounded linear operator on $H$. Suppose that there exists a finite or countably infinite sequence $(F_j)_{j \in I}$ of maps $F_j : \T \to H$ such that $T F_j(\lambda) = \lambda F_j(\lambda)$ for every $j \in I$ and $\lambda \in \T$, and that for every Borel subset $A \subseteq \T$ such that $m(A) = 1$, the subspace $\textrm{span}\{ F_j (\lambda) : j \in I, \lambda \in A \}$ is dense in $H$.
Let $M \sim [0, q, \nu]$ be a complex infinitely divisible random measure on $I \times \Z$ satisfying $\E(M(A)) = 0$ and $\E(\lvert M(A) \rvert^p) \le C' \#(A)$ for every finite subset $A$ of $I \times \Z$ and for some constant $C'>0$ and some $p \in [1,2]$, where $\#$ is the counting measure on $I \times \Z$.
Suppose also that for each $j \in I$, the map $F_j$ is $\beta_j$-Hölderian for some $\beta_j > 1/p - 1/2$. Then the infinitely-divisible distribution of 
\begin{align*}
    \sum_{(j,n) \in I \times \Z} a_j \hat{F_j}(n) M(\{(j,n) \})
\end{align*}
is $T$-invariant and strongly mixing, where $(a_j)_{j \in I}$ is a sequence of positive real numbers such that $\displaystyle \sum_{(j,n) \in I \times \Z} a_j^p \lVert \hat{F_j}(n) \rVert^p < \infty$.
\eco

\medskip
The bound $2/3 < \alpha \le 2$ in Corollary \ref{cormesuresstablesinvholdercashilbert} comes from the inequalities $0 \le 1/\alpha -1/2 <1 $. We now explain how to improve this bound a bit. The idea is to extend the definition of the space $W^{s,2}(\T;H)$ to other values of $s$, following Hörmander \cite[p. 241]{Hormander90}. For $s \in (1,2)$, we denote by $W^{s,2}(\T;H)$ the space of functions $f \in L^2(\T;H)$ such that the weak derivative $f'$ exists and $f' \in W^{s-1,2}(\T;H)$. Let us recall that, for a function $f \in L^2(\T;H)$, the weak derivative $f'$ is the function of $L^2(\T;H)$ satisfying $\int_{\T} f' \varphi \, dm = - \int_{\T} f \varphi' \, dm$ for every infinitely differentiable function $\varphi : \T \to H$. In particular, we have that $\hat{f'}(n) = i n \hat{f}(n)$ for every $n \in \Z$. Let us notice that $f' $ belongs to $W^{s-1,2}(\T;H)$ if and only if $\displaystyle \sum_{n \in \Z} (1+ n ^2)^{s-1} \lVert \hat{f'}(n) \rVert^2 < \infty$, which is equivalent to have $\displaystyle \sum_{n \in \Z} (1+ n ^2)^{s} \lVert \hat{f}(n) \rVert^2 < \infty$, since $(1 + n^2)^{s} \sim (1+ n ^2)^{s-1} n^2$ as $\lvert n \rvert \to \infty$. Hence, if $s \in (1,2)$ satisfies $s > 1/\alpha - 1/2$, i.e., if $s-1 > 1/\alpha -3/2$, and if $f' \in W^{s-1,2}(\T;H)$, then $f \in \mathcal{H}^{s}(\T;H)$ and $\displaystyle \sum_{n \in \Z} \lVert \hat{f}(n) \rVert^\alpha < \infty$. This shows that if $\alpha \in (2/5, 2/3]$, if the $\T$-eigenvector fields $F_j $ in Corollary \ref{cormesuresstablesinvholdercashilbert} are such that the weak derivatives $F_j'$ exist in $L^2(\T;H)$ and are  $\beta_j$-Hölderian for some $\beta_j > 1/\alpha - 3/2$, then the operator $T$ admits a strongly mixing isotropic $\alpha$-stable invariant measure with full support. 

\subsection{Sufficient conditions in a Banach space}
We now consider an infinite-dimensional separable complex Banach space $E$. Of course, we now have to suppose that the Banach space $E$ satisfies some geometric assumptions. In order to control the Fourier coefficients of a function $f \in L^1(\T;E)$, we introduce the Fourier type of a Banach space, following \cite{Konig91}.

\smallskip
The Banach space $E$ is said to be of Fourier type $q \in [1,2]$ if there exists some constant $C>0$ (depending on $q$) such that
\begin{align*}
    (\sum_{n \in \Z} \lVert \hat{f}(n) \rVert^{q'})^{1/q'} \leq C \lVert f \rVert_{L^q(\T;E)}
\end{align*}
for every function $f \in L^q(\T;E)$, where $q'$ is the conjugate exponent of $q$. If $E$ is of Fourier type $q \in [1,2]$, then $E$ is of type $q$ and thus of stable type $\alpha$ for every $0<\alpha<q$. The aim of this subsection is to prove the following result regarding stable invariant measures.

\smallskip

\bth \label{thmchampsvpsunimodstableFouriertype}
Let $E$ be an infinite-dimensional separable complex Banach space of Fourier type $q \in [1,2]$ and of stable type $\alpha \in [1,2]$. Let $T$ be a bounded linear operator on $E$. Suppose that there exists a finite or countably infinite sequence $(F_j)_{j \in I}$ of maps $F_j : \T \to E$ such that $T F_j(\lambda) = \lambda F_j(\lambda)$ for every $j \in I$ and $\lambda \in \T$, and that for every Borel subset $A \subseteq \T$ such that $m(A) = 1$, the subspace $\textrm{span}\{ F_j (\lambda) : j \in I, \lambda \in A \}$ is dense in $E$. Suppose also that for each $j \in I$, the map $F_j$ is $\beta_j$-Hölderian for some $\beta_j > 1/\alpha - 1/q'$, where $q'$ is the conjugate exponent of $q$. Then the operator $T$ admits a strongly mixing isotropic $\alpha$-stable invariant measure with full support.
\eth

\smallskip

Our proof is based on the same ingredients as \cite[Section 3.4]{BayartGrivaux07} and \cite[Thm. 4]{Konig91}. We recall the definition of the Besov space $B_{u,v}^{\lambda}(\T;E)$. For $\lambda \in (0,1)$ and $u,v \geq 1$, the Besov space $B_{u,v}^{\lambda}(\T;E)$ is the space of functions $f \in L^u(\T;E)$ such that
\begin{align*}
    \int_0^{2\pi} \left(\sup_{\lVert y \rVert \leq t} \int_{\T} \lVert f(x+y) - f(x) \rVert^u dx \right)^{v/u} t^{-\lambda v} \frac{dt}{t} < \infty. 
\end{align*}
In the case where $E$ is of Fourier type $q \in [1,2]$, if $\lambda \in (0,1)$, $u,v \geq 1$, and if we define $t$ by the relation
\begin{align*}
    \frac{1}{t} = \lambda + \frac{1}{\max(u', q')},
\end{align*}
then every function $f$ in $B_{u,v}^{\lambda}(\T;E)$ satisfies $\displaystyle \sum_{n \in \Z} \lVert \hat{f}(n) \rVert^v \lvert n \rvert^{v/t -1} < \infty$. Note that, for $\lambda \in (0,1)$ and $u,v \geq 1$, any $\beta$-Hölderian function belongs to $B_{u,v}^\lambda(\T;E)$ when $\beta>\lambda$, by \cite[Proof of Thm. 3.7]{BayartGrivaux07}. 

\smallskip

\bpf[Proof of Theorem \ref{thmchampsvpsunimodstableFouriertype} ]
Suppose that $f : \T \to E$ is a $\beta$-Hölderian function for some $\beta > 1/\alpha - 1/q'$. Taking $\lambda$ such that $\beta > \lambda > 1/\alpha - 1/q'$, as well as $u = 2$ and $v = \alpha$, we obtain that $f$ belongs to $B_{2,\alpha}^{\lambda}(\T;E)$. Defining $t>0$ by $\frac{1}{t} = \lambda + \frac{1}{\max(q',2)} = \lambda + \frac{1}{q'}$, we obtain that $\displaystyle \sum_{n \in \Z} \lVert \hat{f}(n) \rVert^\alpha  \lvert n \rvert^{\alpha/t -1} < \infty$. Since $\alpha/t - 1 = \alpha(\lambda + 1/q') -1 > 1-1 = 0$, we also have that $\displaystyle \sum_{n \in \Z} \lVert \hat{f}(n) \rVert^\alpha < \infty$. The proof of Theorem \ref{thmchampsvpsunimodstableFouriertype} thus follows from Theorem \ref{theoremstablechampsvpsunimod}.
\epf

\smallskip
Similarly, for infinitely divisible invariant measures, we have the following result. In particular, this recovers the result of Bayart and Grivaux \cite[Thm. 3.7]{BayartGrivaux07}.

\smallskip

\bth \label{LoiIDinvarianteFouriertype}
Let $E$ be an infinite-dimensional separable complex Banach space of Fourier type $r \in [1,2]$ and of type $p \in [1,2]$. Let $T$ be a bounded linear operator on $E$. Suppose that there exists a finite or countably infinite sequence $(F_j)_{j \in I}$ of maps $F_j : \T \to E$ such that $T F_j(\lambda) = \lambda F_j(\lambda)$ for every $j \in I$ and $\lambda \in \T$, and that for every Borel subset $A \subseteq \T$ such that $m(A) = 1$, the subspace $\textrm{span}\{ F_j (\lambda) : j \in I, \lambda \in A \}$ is dense in $E$. Suppose also that for each $j \in I$, the map $F_j$ is $\beta_j$-Hölderian for some $\beta_j > 1/p - 1/r'$, where $r'$ is the conjugate exponent of $r$.
Let $M \sim [0, q, \nu]$ be a complex infinitely divisible random measure on $I \times \Z$ satisfying $\E(M(A)) = 0$ and $\E(\lvert M(A) \rvert^p) \le C' \#(A)$ for every finite subset $A$ of $I \times \Z$ and for some constant $C'>0$, where $\#$ is the counting measure on $I \times \Z$.
Then the infinitely-divisible distribution of 
\begin{align*}
    \sum_{(j,n) \in I \times \Z} a_j \hat{F_j}(n) M(\{(j,n) \})
\end{align*}
is $T$-invariant and strongly mixing, where $(a_j)_{j \in I}$ is a sequence of positive real numbers such that $\displaystyle \sum_{(j,n) \in I \times \Z} a_j^p \lVert \hat{F_j}(n) \rVert^p < \infty$. In particular, $T$ admits a strongly mixing isotropic Gaussian invariant measure with full support.
\eth

\subsection{Examples}
We now study some examples of operators on Banach spaces admitting isotropic stable or infinitely divisible ergodic invariant measures. We first start with weighted shifts on the complex spaces $\ell_p(\Z_+)$ and $c_0(\Z_+)$ when $1 \leq p < \infty$. If $w=(w_n)_{n \geq 1}$ is a bounded sequence of nonzero complex numbers, the weighted shift on these spaces is the operator defined by $B_w e_0 = 0$ and $B_w e_n = w_n e_{n-1}$ for $n \geq 1$, where $(e_n)_{n \geq 0}$ denotes the canonical basis of these spaces.

\smallskip

\bpr \label{exemplestableshiftsàpoids}
Let $E = \ell_p(\Z_+)$ or $E = c_0(\Z_+)$ for $1 \leq p < \infty$, and let $w = (w_n)_{n \geq 1} $ be a bounded sequence of nonzero complex numbers. Suppose that $0 <\alpha < p $ if $p < 2$ and that $0 < \alpha \leq 2 $ if $p \geq 2$ in the case where $E = \ell_p(\Z_+)$, and that $0< \alpha < 1$ in the case where $E = c_0(\Z_+)$. 
If 
\begin{align*}
    \sum_{n \geq 1} \frac{1}{\lvert w_1 \dotsm w_n \rvert^\alpha} < \infty,
\end{align*}
then $B_w$ admits a strongly mixing isotropic $\alpha$-stable invariant measure with full support. 
\epr

\smallskip

\bpf
We prove Proposition \ref{exemplestableshiftsàpoids} for $E = \ell_p(\Z_+)$, since the proof is essentially the same for $c_0(\Z_+)$.
Let $\beta_0 = 1$ and $\beta_n = w_1 \dotsm w_n$ if $n \geq 1$. 
The following $\T$-eigenvector field of $B_w$ given by 
\begin{align*}
    F : \lambda \mapsto \displaystyle \sum_{n \geq 0} \frac{\lambda^n}{\beta_n} e_n
\end{align*}
is well-defined since $\ell_\alpha(\Z_+) \subseteq \ell_p(\Z_+)$ if $0 < \alpha \leq p $. In this case, the Fourier coefficients of $F$ are given by $\hat{F}(n) = 0$ if $n < 0$ and by $\hat{F}(n) = \frac{e_n}{\beta_n}$ if $n \geq 0$. Indeed, for each $y^* \in (\ell_p(\Z_+))^* = \ell_{p'}(\Z_+)$, we have
\begin{align*}
    \lan y^*, \hat{F}(n) \ran &= \int_{\T} \overline{\lambda}^n \sum_{k \geq 0} \frac{\lambda^k}{\beta_k} \lan y^*, e_k \ran dm(\lambda) \\
    &= \sum_{k \geq 0} \frac{y_k}{\beta_k} \delta_{k,n} 
\end{align*}
since the series $\displaystyle \sum_{n \geq 0} \frac{\lambda^n}{\beta_n} y_n$ is uniformly convergent on $\T$ by Hölder's inequality. The map $F$ obviously belongs to $L^1(\T;E)$ and the end of the proof follows from Theorem \ref{theoremstablechampsvpsunimod} and from the fact that $F$ is spanning with respect to the normalized Lebesgue measure on $\T$, as well as the fact that $E$ is of stable type $\alpha $ if $p<2$ and of stable type $2$ if $p \geq 2$.
\epf

\smallskip
Let us note that it is possible to formulate an analogue of Proposition \ref{exemplestableshiftsàpoids} for infinitely divisible distributions, by supposing that $\sum_{n \geq 1} \frac{1}{\lvert w_1 \dotsm w_n \rvert^2} < \infty$ if $p \ge 2$ and $\sum_{n \geq 1} \frac{1}{\lvert w_1 \dotsm w_n \rvert^p} < \infty$ if $1 \le p <2 $. However, for Gaussian measures, we have the following improvement, which recovers the classical criterion for weighted shifts \cite[Example 5.5.2]{BM09}.
\smallskip

\bpr \label{mesurestableinvshiftspoids}
Let $E = \ell_p(\Z_+)$ or $E = c_0(\Z_+)$ for $1 \leq p < \infty$, and let $w = (w_n)_{n \geq 1} $ be a bounded sequence of nonzero complex numbers. 
Suppose that
\begin{align*}
    \sum_{n \geq 1} \frac{1}{\lvert w_1 \dotsm w_n \rvert^p} < \infty,
\end{align*}
in the case where $E = \ell_p(\Z_+)$, and that
\begin{align*}
    \sum_{n \geq 1} \frac{1}{\lvert w_1 \dotsm w_n \rvert} < \infty,
\end{align*}
in the case where $E = c_0(\Z_+)$.
Then $B_w$ admits a strongly mixing isotropic Gaussian invariant measure with full support. 
\epr

\smallskip

\bpf
We again use the expression of the Fourier coefficients of the $\T$-eigenvector field $F$ obtained in the proof of Proposition \ref{exemplestableshiftsàpoids}. We obtained that $\hat{F}(n) = 0$ if $n <0$ and $\hat{F}(n) = \frac{e_n}{\beta_n}$ if $n \ge 0$. Suppose that $E = \ell_p(\Z_+)$. Then, for $p \in [1,2]$, the Banach space $E$ is of type $p$ and $\displaystyle \sum_{n \in \Z} \lVert \hat{F}(n) \rVert^p < \infty$ by the assumption on the weight $(w_n)_{n \ge 1}$. Hence, by Corollary \ref{corexistencemesureGaussinvtypep}, the operator $B_w$ admits a strongly mixing isotropic Gaussian invariant measure with full support. It remains to show it for $p>2$. This time, the Banach space $E$ is of type $2$. In order to make use of the summability hypothesis on the weight $(w_n)_{n \ge 0}$, we must redefine the random integral $\int_{\Z_+} \hat{F}(n) dM(n)$, where $M$ is a complex isotropic Gaussian random measure on $\Z_+$, based on \cite[p. 121]{BM09}. For $N \ge 0$, let us set $f_N$ by $f_N(n)  = \hat{F}(n)$ if $0 \le n \le N$ and $f_N(n) = 0 $ if $n >N$. Then, for every integers $N > M$, we have that
\begin{align*}
    \E  \left \lVert \sum_{n=M+1}^N \hat{F}(n) M(\{n \}) \right \rVert^p = \E(\lvert M(\{0 \}) \rvert^p) \sum_{n=M+1}^N \frac{1}{\lvert w_1 \dotsm w_n \rvert^p}
\end{align*}
which converges to $0$ as $M,N \to \infty$. Hence, the sequence of the random integrals $\int_{\Z_+} f_N \, dM$ converges in $L^p(\Omg; E)$ to $\int_{\Z_+} \hat{F}(n) dM(n)$. Moreover, using the convergence of $(f_N)_{N \ge 0}$ to $\hat{F}$ in $L^p(\Z_+;E)$, the characteristic function of $\int_{\Z_+} \hat{F}(n) dM(n)$ is given for every $x^* \in E^*$ by
\begin{align*}
    \E(e^{i \Re(\lan x^*, \int_{\Z_+} \hat{F}(n) dM(n) \ran)}) = \exp(-\frac{1}{2} \sum_{n \ge 0} \lvert \lan x^*, \hat{F}(n) \ran \rvert^2),
\end{align*}
where the series is convergent by Hölder's inequality and by the summability condition on the weight $(w_n)_{n \ge 1}$. This shows in this case that the isotropic Gaussian distribution of $\int_{\Z_+} \hat{F}(n) dM(n)$ is invariant and strongly mixing with respect to the operator $B_w$ with full support. This concludes the proof in the case where $E = \ell_p(\Z_+)$. The case where $E = c_0(\Z_+)$ can be treated similarly.
\epf

\smallskip
We now study adjoints of multiplication operators on the Hardy space $H^2(\D)$. Let us recall that if $\phi \in H^\infty(\D)$, the multiplication operator by $\phi$ is given by $M_\phi : f \mapsto f \phi$ on $H^2(\D)$. 

\smallskip

\bpr \label{exempleadjointopmultiplication}
Let $E = H^2(\D)$ and let $\phi \in H^\infty(\D)$ nonconstant such that $\phi(\D) \cap \T \ne \emptyset$. Then for every $\alpha \in (0,2]$, the adjoint $M_{\phi}^*$ admits a strongly mixing isotropic $\alpha$-stable invariant measure with full support.  
\epr

\smallskip

\bpf
We mix the proof given in \cite[Example 5.5.1]{BM09} and \cite[Proof of Thm. 8.2]{BayartGrivauxMatheronMenet25}.
Let $a \in \D$ such that $\phi(a) \in \T$. Let $W$ be an open neighbor of $a$ in $\D$. Then $\phi(W)$ is an open set intersecting $\T$, so $\phi(W) \cap \T$ is uncountable. Thus, there exists $b \in W$ such that $\phi(b) \in \T$ and $\phi'(b) \ne 0$. By the inverse function theorem, there exists an open set $U$ in $\D$ with $b \in U$ such that $\phi : U \to \phi(U)$ is a diffeomorphism. Moreover, taking $U $ smaller if necessary, we may assume that $I := \phi(U) \cap \T$ is a proper open arc of $\T$. Let $\psi$ be the inverse of $\phi : U \to \phi(U)$. Note that $\psi$ is holomorphic. Let us also denote by $J := \{ \overline{z} : z \in I \}$. Consider $\chi \in \mathcal{C}^\infty(\T)$ a \enquote{cut-off} function such that $\chi = 0$ outside $J$ and $\chi >0$ on $J$. Then, the map $F : \T \to H^2(\D)$ defined by $F(e^{i \theta}) = \chi(e^{i \theta}) k_{\psi(e^{-i \theta})}$ if $e^{i \theta} \in J$ and $F(e^{i \theta})=0$ if $e^{i \theta} \notin J$ belongs to $\mathcal{C}^{\infty}(\T;H^2(\D)$ and it satisfies $M_{\phi}^* E(\lambda) = \lambda E(\lambda) $ for every $\lambda \in \T$. In particular, we have that $\lVert \hat{F}(n) \rVert = O (1/\lvert n \rvert^k)$ for every integer $k \geq 1$ and thus $\displaystyle \sum_{n \in \Z} \lVert \hat{F}(n) \rVert^\alpha < \infty$. Moreover, the $\T$-eigenvector field $F$ is spanning with respect to the normalized Lebesgue measure on $\T$. This concludes the proof of Proposition \ref{exempleadjointopmultiplication}.
\epf
\smallskip
Similarly, we can formulate an analogue of Proposition \ref{exempleadjointopmultiplication} for infinitely disivible distributions, since the $\T$-eigenvector field is in $\mathcal{C}^{\infty}(\T; H^2(\D))$.

\smallskip
We now study Kalisch-type operators on $L^p[0,2\pi]$.

\smallskip

\bpr \label{Kalisch-typestableLP}
Let $E = L^p[0,2\pi]$ for $1 < p < \infty$, and let $T$ be the operator on $E$ defined by
\begin{align*}
    Tf(\theta) = e^{i \theta} f(\theta) - \int_0^{\theta} i e^{it} f(t) dt.
\end{align*}
Let $\alpha \in (1,2]$ if $p \geq 2$ and $\alpha \in (1,p)$ if $1 <p <2$. Then, $T$ admits a strongly mixing isotropic $\alpha$-stable invariant measure with full support. 

Finally, for every $1 \leq p < \infty$, $T$ admits a strongly mixing isotropic Gaussian invariant measure with full support.
\epr

\smallskip

\bpf
It is well-known (see, e.g., \cite[Example 5.5.3]{BM09}) that for every $u \in [0, 2\pi)$, the function $F(e^{iu}) := \indic_{(u, 2\pi)}$ satisfies $T F(e^{iu}) = e^{iu} F(e^{iu})$. Moreover, it is well-known that $F$ is $1/p$-Hölderian. Note that $E$ has type $\min(p,2)$ (and thus stable type $\alpha$ for every $\alpha \in (0,p)$ if $p<2$ and $\alpha \in (0,2]$ if $p \ge 2$), and Fourier-type $p$ if $p \in [1,2]$. Thus, the inequality $1/p > 1/\alpha - 1/p'$ is satisfied if $\alpha > 1$ and we conclude by Theorem \ref{thmchampsvpsunimodstableFouriertype}  in the case where $p \in (1,2]$ that $T$ admits a strongly mixing isotropic $\alpha$-stable invariant measure with full support. Note that it only covers the case where $p \in (1,2]$. Fortunately, it is possible here to compute the Fourier coefficients of $F$. 

For $n \in \Z$, let $e_n : x \mapsto e^{-inx}$ on $[0, 2\pi]$. Then 
\[ \hat{F}(n) = \frac{1-e_n}{2\pi in} \quad \text{if } \lvert n \rvert \ge 1 \quad \text{and} \quad \hat{F}(0) : x \mapsto \frac{1}{2\pi} x. \]
Since the $\T$-eigenvector field $F$ is spanning with respect to the normalized Lebesgue measure on $\T$, the subspace $\textrm{span}\{ u_n : n \in \Z \}$ is dense in $E$, where $u_n := \hat{F}(n)$ for every $n \in \Z$. Moreover, for $n \ne 0$, we have
\begin{align*}
    \lVert u_n \rVert^p &= \int_0^{2\pi} \frac{1}{(2\pi)^p} \frac{\lvert e^{-inx} - 1 \rvert^p}{\lvert n \rvert^p} dx \\
    &=\frac{2^p}{(2\pi)^p \lvert n \rvert^p} \int_0^{2 \pi} \left \lvert \sin(\frac{nx}{2}) \right \rvert^p dx \\
    &= \frac{2}{\lvert n \rvert^p \pi^p} \int_{0}^{\lvert n \rvert \pi} \lvert \sin(y) \rvert^p \frac{dy}{\lvert n \rvert} \\
    &= \frac{2 C_p}{\pi^p \lvert n \rvert^p},
\end{align*}
where $C_p $ is a positive constant depending on $p$ but not on $n$.
Thus, if $\alpha >1$, then the series $\displaystyle \sum_{n \ne 0} \lVert u_n \rVert^\alpha$ is convergent and we conclude that the operator $T$ admits a strongly mixing isotropic $\alpha$-stable invariant measure when $\alpha \in (1,2]$ and $p > 2$. 

\smallskip
Finally, suppose that $1 \le p < \infty$ and let us note that $E$ has Fourier type $\min(p,p')$ and type $\min(2,p)$. Hence, if $p \in [1,2]$, then the inequality $1/p > 1/p - 1/p'$ is satisfied. Thus, by Theorem \ref{LoiIDinvarianteFouriertype}, the operator $T$ admits a strongly mixing isotropic Gaussian invariant measure with full support. This concludes the proof of Proposition \ref{Kalisch-typestableLP}.
\epf

Let us note that the operator introduced in Proposition \ref{Kalisch-typestableLP} has no $\T$-eigenvector in the space $\mathcal{C}_0[0,2\pi]$ of complex continuous functions on $[0,2\pi]$ vanishing at $0$. We now show that, nevertheless, it admits a strongly mixing $\alpha$-stable isotropic invariant measure when $\alpha >1$.

\smallskip

\bpr \label{Kalisch-typestableC0}
Let $E = \mathcal{C}_0[0,2\pi]$ and let $T$ be the operator on $E$ defined by
\begin{align*}
    Tf(\theta) = e^{i \theta} f(\theta) - \int_0^{\theta} i e^{it} f(t) dt.
\end{align*}
Let $\alpha \in (1,2]$. Then, $T$ admits a strongly mixing isotropic $\alpha$-stable invariant measure with full support.
\epr

\smallskip

\bpf
Let $u_n = \frac{1-e_n}{2 i \pi n}$ if $\lvert n \rvert \ge 1$ and $u_0 : x \mapsto  \frac{1}{2\pi} x$, where again $e_n : x \mapsto e^{-inx}$.
The $u_n$'s are de Fourier coefficients of the $\T$-eigenvector field of $T$ computed in the proof of Proposition \ref{Kalisch-typestableLP}. Even though the $\T$-eigenvector field does not belong to $E$, its Fourier coefficients $(u_n)_{n \in \Z}$ belong to $E$.

\smallskip
Let $(\theta_k^{(\alpha)})_{k \in \Z}$ be a sequence of i.i.d real standard $\alpha$-stable random variables.
Then from Marcus and Pisier \cite[p. 294]{MarcusPisier84}, the series $\displaystyle \sum_{k \ge 1} \frac{\theta_{-k}^{(\alpha)}}{2i\pi k} e^{ikt}$ is almost surely uniformly convergent on $[0,2\pi]$ when $\alpha>1$, and moreover, the series $\displaystyle \sum_{k \ge 1} \frac{\theta_{-k}^{(\alpha)}}{2 i \pi k}$ is almost surely convergent too when $\alpha>1$, since $\C$ has stable type $2$. Indeed, the criterion given in \cite[p. 294]{MarcusPisier84} states that $\displaystyle \sum_{k \ge 1} \frac{\theta_{-k}^{(\alpha)}}{2i\pi k} e^{ikt}$ is almost surely uniformly convergent on $[0,2\pi]$ if and only if $\displaystyle \sum_{n \ge 1} \frac{(\sum_{k \ge n} \lvert \frac{1}{2 i \pi k} \rvert^\alpha)^{1/\alpha}}{n (\log(n+1))^{1/\alpha}} < \infty$ for $\alpha \in [1,2]$, and this condition is satisfied when $\alpha >1$. Since $\displaystyle \sum_{k \ge 1} \frac{\theta_{k}^{(\alpha)} e^{-ikt}}{2 i \pi (-k)} = \overline{\sum_{k \ge 1} \frac{\theta_{k}^{(\alpha)} e^{ikt}}{2 i \pi k}}$ and $\displaystyle \sum_{k \ge 1} \frac{\theta_{k}^{(\alpha)}
}{2 i \pi (-k)} = \overline{\sum_{k \ge 1} \frac{\theta_{k}^{(\alpha)}
}{2 i \pi k}}$, these two series are also almost surely uniformly convergent on $[0,2\pi]$. Let $(X_n)_{n \in \Z}$ be a sequence of i.i.d complex isotropic $\alpha$-stable random variables with characteristic function $\E(e^{i \Re(\overline{z} X_n)}) = e^{-\lvert z \rvert^\alpha}$. Then, almost surely, the series $\displaystyle \sum_{n \in \Z} X_n u_n$ is convergent in $E$ when $\alpha > 1$. The distribution of this random variable is easily seen to be $T$-invariant, and also strongly mixing by the proof of Theorem \ref{theoremstablechampsvpsunimod}. Let us now compute the characteristic function of the random vector $G := \displaystyle \sum_{n \in \Z} X_n u_n$.
Let us note that $\displaystyle \sum_{\lvert n \rvert \le N} X_n u_n$ converges in distribution to $G$ as $N \to \infty$. Hence, for $x^* \in E^*$, we have that
\begin{align*}
    \E[e^{i \Re(\lan x^*, G \ran)}] &= \lim_{N \to \infty} \E[e^{i \Re(\lan x^*, \sum_{\lvert n \rvert \le N} X_n u_n \ran)} ] \\
    &=\lim_{N \to \infty} \prod_{\lvert n \rvert \le N} \E[e^{i \Re(\lan x^*, u_n \ran X_n)}] \\
    &= \lim_{N \to \infty} \exp(-\sum_{\lvert n \rvert \le N} \lvert \lan x^*, u_n \ran \rvert^\alpha) \\
    &= \exp(-\sum_{n \in \Z} \lvert \lan x^*, u_n \ran \rvert^\alpha),
\end{align*}
where we used the fact that since $\alpha >1$, the series $\displaystyle \sum_{n \in \Z} \lvert \lan x^*, u_n \ran \rvert^\alpha$ is convergent. This shows that the distribution of $G$ is isotropic $\alpha$-stable. In order to show that $G$ has full support, it remains to show that the subspace $\textrm{span}\{ u_n : n \in \Z \}$ is dense in $E$.

\smallskip
Let $L $ be a linear functional on $E$ such that $L = 0$ on $\textrm{span}\{ u_n : n \in \Z \}$. By the Hahn-Banach theorem, $L$ extends to a linear functional on $\mathcal{C}[0,2\pi]$ and thus there exists a complex measure on $[0,2\pi]$ such that $L(f) = \int_{[0,2\pi]} f d\mu$ for every function $f \in \mathcal{C}[0,2\pi]$. Let $\nu$ be the restriction of $\mu$ to $(0,2\pi]$. Then, for every function $f \in E$, we have that $L(f) = \int_{(0,2\pi]} f d\nu$. Since $L(u_n) = 0$ for every $n \ne 0$, we obtain that $\hat{\nu}(n) = \nu(0,2\pi]$ for every $n \in \Z$, where $\hat{\nu}(n)$ is the $n$-th Fourier coefficient of the complex measure $\nu$. Let $\beta := \nu(0,2\pi]$. Since the Fourier coefficients of the complex measure $\nu - \beta \delta_{2\pi}$ are all equal to zero, we obtain that $\nu = \beta \delta_{2\pi}$. Finally, since $L(u_0) = 0$, we obtain that $\beta = 0$. This shows that $L = 0$, and that the subspace $\textrm{span}\{ u_n : n \in \Z \}$ is dense in $E$. This concludes the proof of Proposition \ref{Kalisch-typestableC0}. 
\epf

\section{Ergodic stable and infinitely divisible invariant measures for strongly continuous semigroups of operators}

We now come to semigroups of operators on $E$ whose generators admit a perfectly spanning family of purely imaginary eigenvector fields. By \cite[Thm. 3.8]{Gill1}, if $E$ is of stable type $\alpha \in (0,2)$, and if $(T_t)_{t \geq 0}$ is a strongly continuous semigroup of operators on $E$ such that there exists an $\alpha$-integrable family $(u_t)_{t \in J}$ of purely imaginary eigenvectors for the generator (with $J$ an interval of $\R$) for which $\textrm{span}\{ u_t : t \in J \setminus B \}$ is dense for every Borel subset $B$ with measure zero, then the distribution of the random vector $\int_J u_t M(dt)$ is a non-ergodic $\alpha$-stable invariant measure with full support for $(T_t)_{t \ge 0}$, for every complex isotropic $\alpha$-stable random measure on $J$. However, in the case where $\alpha = 2$, this distribution is in fact strongly mixing with respect to the semigroup by \cite[Thm. 4.1]{AGNEESSENS_2026}, \cite[Thm. 3.2]{ChakirMourchid18} and \cite[Thm. 3.8]{Gill1}. In the case of operators, we made use of the Fourier coefficients of the $\T$-eigenvector fields. For semigroups, it is natural to make use of the Fourier transforms of the purely imaginary eigenvector fields of the generators. The idea to make use of the Fourier transform was also used in \cite[Thm. 4.2, p. 1139]{AGNEESSENS_2026}. If $F : \R \to E$ is in $L^1(\R;E)$, the Fourier transform of $F$ is the map $\hat{F} :  x \in \R \mapsto \int_{\R} F(t) e^{-ixt} dt$.  

\subsection{General results}
The main results for the existence of ergodic stable or infinitely divisible invariant measures for semigroups of operators whose generators admit a sequence of perfectly spanning purely imaginary eigenvector fields are the following. First, we focus on stable distributions.
Let us recall that if $(T_t)_{t \geq 0}$ is a $\mathcal{C}_0$-semigroup of operators on a Banach space $E$ with generator $A$, then $A x_\lambda =  \lambda x_\lambda$ with $\lambda \in \C$ and $x_\lambda \in E$ implies that $T_t x_\lambda = e^{\lambda t} x_\lambda$ for any $t \ge 0$.

\medskip
\bth \label{thminvmeasurestablesemigroups}
Let $E$ be an infinite-dimensional separable complex Banach space of stable type $\alpha \in (0,2]$ and let $(T_t)_{t \geq 0}$ be a $\mathcal{C}_0$-semigroup of operators on $E$. Let $A$ be the generator of this semigroup. Suppose that there exists a finite or countably infinite sequence $(F_j)_{j \in I}$ of maps in $L^1(J;E)$ for some interval $J \subseteq \R$ such that $A F_j(t) = i t F_j(t)$ for every $j \in I$ and $t \in J$, and that the subspace $\textrm{span}\{ F_j(t) : j \in I, \, t \in J \setminus B \}$ is dense in $E$ for every Borel subset $B$ of $J$ with Lebesgue measure zero. Suppose also that $ \int_{\R} \lVert \hat{F_j}(t) \rVert^\alpha dt < \infty$ for every $j \in I$, where $\hat{F_j}$ denotes the Fourier transform of the map $F_j$. Then the semigroup $(T_t)_{t \geq 0}$ admits a strongly mixing isotropic $\alpha$-stable invariant measure with full support.
\eth

\medskip

\bpf
Let us note that the maps $t \mapsto \hat{F_j}(t)$ are well-defined since the maps $F_j$ are integrable on $J$. Moreover, the subspace $\textrm{span}\{ \hat{F_j}(t) : j \in I, t \in \R \setminus B \}$ is dense in $E$ for every Borel subset $B$ of $\R$ with Lebesgue measure zero. Indeed, let $x^* \in E^*$ vanishing on this subspace. Then
\begin{align*}
    \int_{\R} \lan x^*, F_j(\theta) \ran \indic_{J}(\theta) e^{-it \theta} d\theta = 0
\end{align*}
for every $j \in I$ and $t \in \R \setminus B$. Since $F_j \in L^1(J;E)$ for every $j \in I$, this implies that for every $j \in I$, there exists a Borel subset $C_j$ of $J$ with Lebesgue measure zero such that $\lan x^*, F_j(\theta) \ran = 0$ for every $\theta \in J \setminus C_j$. Setting $C := \displaystyle  \bigcup_{j \in I} C_j$, which has Lebesgue measure zero in $J$, we obtain that $\lan x^*, F_j(\theta) \ran = 0$ for every $\theta \in J \setminus C$ and for every $j \in I$.
Thus $x^*$ vanishes on the dense subspace $\textrm{span}\{ F_j(t) : j \in I, \, t \in J \setminus C \}$ and $x^* = 0$. 

\smallskip
Next, let us notice that, for every $s \geq 0$, $j \in I$ and $t \in \R$, we have that
\begin{align*}
    T_s \hat{F_j}(t) = \int_J e^{-it \theta} T_s F_j(\theta) d\theta = \int_J e^{-it \theta}  e^{i s \theta} F_j(\theta) d\theta = \hat{F_j}(t-s).
\end{align*}

\medskip
Since $\int_{\R} \lVert \hat{F_j}(t) \rVert^\alpha dt < \infty$ for every $j \in I$, we can choose a sequence $(a_j)_{j \in I}$ of positive real numbers such that $\displaystyle \sum_{j \in I} a_j^\alpha \int_{\R} \lVert \hat{F_j}(t) \rVert^\alpha dt < \infty$. Let $M$ be a complex isotropic $\alpha$-stable random measure on the $\sigma$-finite space $(I \times \R, \mathcal{P}(I) \otimes \mathcal{\textrm{Bor}}(\R), \# \otimes m)$, where $\#$ is the counting measure on $I$ and $m$ is the Lebesgue measure on $\R$. Hence, the random vector
\begin{align*}
    G := \int_{I \times \R} a_j \hat{F_j}(t)\, dM(j,t)
\end{align*}
is well-defined, and denoting by $\mu_G$ the distribution of $G$, its characteristic function is given for every $x^* \in E^*$ by
\begin{align} \label{fonccaractrandomvectorsemigroups}
    \hat{\mu_G}(x^*) = \exp \left(-c_0 \sum_{j \in I} a_j^\alpha \int_{\R} \lvert \lan x^*, \hat{F_j}(t) \ran \rvert^\alpha dt \right).
\end{align}
Moreover, $G$ is an isotropic $\alpha$-stable random vector, and using its characteristic function (\ref{fonccaractrandomvectorsemigroups}), together with the fact that $\textrm{span}\{ \hat{F_j}(t) : j \in I, t \in \R \setminus B \}$ is dense in $E$ for every Borel subset $B$ of $\R$ with Lebesgue measure zero, we can show that $G$ has full support. It is also easy to see that $\mu_G$ is invariant under $(T_s)_{s \ge 0}$ using (\ref{fonccaractrandomvectorsemigroups}) and the fact that $T_s \hat{F_j}(t) = \hat{F_j}(t-s)$ for every $j \in I$, $s \ge 0$ and $t \in \R$.

\smallskip
It only remain to show that $\mu_G$ is strongly mixing with respect to $(T_s)_{s \ge 0}$. To do so, we just have to notice that the proof of mixing given in \cite[Thm. 3.1]{Gill1} can be adapted to the product space $I \times \R$. Indeed, notice that, for $N >0$ and for $0<q< \alpha$, we have by (\ref{inequalityqmomentstable}) that 
\begin{align*}
    \left \lVert \int_{I \times [N, \infty)} a_j \hat{F_j}(t) \, dM(j,t) \right \rVert_q^{\alpha} &\leq C \int_{I \times [N,\infty)} a_j^{\alpha} \lVert \hat{F_j}(t) \rVert^\alpha d(\# \otimes m)(j,t) \\ &\le C \int_{N}^\infty \sum_{j \in I} a_j^{\alpha} \lVert \hat{F_j}(t) \rVert^\alpha dt
\end{align*}
for some constant $C$ depending on $q$ and $\alpha$ but not on $N$, so that $\int_{I \times [0, N)} a_j \hat{F_j}(t) \, dM(j,t)$ converges in probability to $\int_{I \times [0, \infty)} a_j \hat{F_j}(t) \, dM(j,t)$ as $N \to \infty$, by Lebesgue's convergence theorem. Moreover, for $N >0$ and $s > 2N$, the random vectors $\int_{I \times (-\infty, -N+s)} a_j \hat{F_j}(t-s) \, dM(j,t)$ and $\int_{I \times (-\infty,-N)} a_j \hat{F_j}(t) \, dM(j,t)$ have the same distribution (just by considering their characteristic functions). And finally, for $N >0$ and $s > 2N$, since $I \times (-\infty, -N+s)$ and $I \times [-N+s, \infty)$ are disjoint, the random integrals $\int_{I \times (-\infty, -N+s) } a_j \hat{F_j}(t-s) \, dM(j,t)$ and $\int_{I \times [-N+s, \infty)} a_j \hat{F_j}(t-s) \, dM(j,t)$ are independent. This is enough to adapt the proof of mixing given in \cite[Thm. 3.1]{Gill1}.
 
\smallskip
This concludes the proof of Theorem \ref{thminvmeasurestablesemigroups}. 
\epf

\smallskip
For infinitely divisible distributions, we have the following result. 

\smallskip

\bth
Let $E$ be an infinite-dimensional separable complex Banach space of type $p \in [1,2]$ and let $(T_t)_{t \geq 0}$ be a $\mathcal{C}_0$-semigroup of operators on $E$. Let $A$ be the generator of this semigroup. Suppose that there exists a finite or countably infinite sequence $(F_j)_{j \in I}$ of maps in $L^1(J;E)$ for some interval $J \subseteq \R$ such that $A F_j(t) = i t F_j(t)$ for every $j \in I$ and $t \in J$, and that the subspace $\textrm{span}\{ F_j(t) : j \in I, \, t \in J \setminus B \}$ is dense in $E$ for every Borel subset $B$ of $J$ with Lebesgue measure zero. Suppose also that $ \int_{\R} \lVert \hat{F_j}(t) \rVert^p dt < \infty$ for every $j \in I$, where $\hat{F_j}$ denotes the Fourier transform of the map $F_j$. Let $M \sim [0, q, \nu]$ be a complex infinitely divisible random measure on $I \times \R$ satisfying $\E(M(A)) = 0$ and $\E(\lvert M(A) \rvert^p) \le C' (\# \otimes m)(A)$ for every measurable subset $A$ of $I \times \R$ such that $(\# \otimes m)(A) <\infty $ and for some constant $C'>0$, where $\#$ is the counting measure on $I$ and $m$ is the Lebesgue measure on $\R$. 
Then the infinitely divisible distribution of the random vector
\begin{align*}
    \int_{I \times \R} a_j \hat{F_j}(t) \, dM(j,t)
\end{align*}
is invariant and strongly mixing with respect to the
semigroup $(T_t)_{t \geq 0}$, where $(a_j)_{j \in I}$ is a sequence of positive real numbers such that $\displaystyle \sum_{j \in I} a_j^p \int_{\R} \lVert \hat{F_j}(t) \rVert^p dt < \infty$.
\eth

\bigskip

We now discuss sufficient conditions in a Hilbert space in order to satisfy Theorem \ref{thminvmeasurestablesemigroups}. For Banach spaces, it is also possible to define the Fourier type exactly as for the operator case. However, for the Fourier transform on $\R$, the relations between the Besov space and the Lorentz space are more complicated than in the discrete case on $\Z$. A reference on the subject is the article \cite[Section 10]{Fouriersurvey98}. For this reason, we won't come into the Banach space setting.

\subsection{Sufficient conditions in a Hilbert space}
Let $H$ be an infinite-dimensional separable complex Hilbert space. To keep this paper concise, we state the results only for stable measures. Analogous results for infinitely divisible measures can be formulated.

For every function function $f \in L^1(\R;H) \cap L^2(\R;H)$, the Plancherel theorem yields that the $L^2$-norm of $f$ is equal to the norm of the Fourier transform of $f$ (up to a multiplicative constant), see \cite[Theorem 2.4.9]{HytonenvanNeerverVeraarWeis16}. This gives the following result regarding Gaussian measures. 

\smallskip

\bpr \label{propexistencemesureinvgaussiennesemigroupe}
Let $H$ be an infinite-dimensional separable complex Hilbert space and let $(T_t)_{t \geq 0}$ be a $\mathcal{C}_0$-semigroup of operators on $H$. Let $A$ be the generator of this semigroup. Suppose that there exists a finite or countably infinite sequence $(F_j)_{j \in I}$ of maps in $L^1(J;H) \cap L^2(J;H)$ for some interval $J \subseteq \R$ such that $A F_j(t) = i t F_j(t)$ for every $j \in I$ and $t \in J$, and that the subspace $\textrm{span}\{ F_j(t) : j \in I, \, t \in J \setminus B \}$ is dense in $H$ for every Borel subset $B$ of $J$ with Lebesgue measure zero. Then the semigroup $(T_t)_{t \geq 0}$ admits a strongly mixing isotropic Gaussian invariant measure with full support.
\epr

\medskip
Similarly to the case of operators, we can introduce Sobolev spaces in order to control the Fourier transforms of the purely imaginary eigenvector fields of the generators. More precisely, following \cite[Section 7.9]{Hormander90}, we denote by $\mathcal{H}^{s}(\R;H)$ the space of all functions $f \in L^1(\R;H) \cap L^2(\R; H)$ such that $$\int_{\R} (1+ \lvert \xi \rvert^2)^s \lVert \hat{f}(\xi)\rVert^2 d\xi < \infty.$$
Similarly to the torus case, one can show using Plancherel's identity that, for $s \in (0,1)$, the space $\mathcal{H}^{s}(\R;H)$ coincides with the space $W^{s,2}(\R;H)$ of functions $f$ in $L^1(\R;H) \cap L^2(\R;H)$ such that $$\int \int_{\R^2} \frac{\lVert f(x) - f(y) \rVert^2}{\lvert x-y \rvert^{1+2s}} dx dy = \int \int_{\R^2} \frac{\lVert f(x+h) - f(x) \rVert^2}{\lvert h \rvert^{1+2s}} dx dh < \infty.$$ 
Moreover, for $\alpha \in (0,2]$, if there exists $s > 1/\alpha -1/2$ such that $f \in \mathcal{H}^{s}(\R;H)$, then $\int_{\R} \lVert \hat{f}(\xi) \rVert^\alpha d\xi < \infty$, by Hölder's inequality. Finally, every $\beta$-Hölderian function $f \in L^1(\R;H) \cap L^2(\R;H)$ belongs to $\mathcal{H}^{s}(\R;H)$ when $0 < 2s <\beta < 1$. 
Indeed, let us notice that, for every $\lvert h \rvert >1$, we have that 
$$\int_{\R} \lVert f(x+h)-f(x) \rVert^2 dx \le 2 \lVert f \rVert_2^2.$$ 
And for every $\lvert h \rvert \le 1$ with $h \ne 0$, we have that 
\begin{align*}
    \int_{\R} \lVert f(x+h) - f(x) \rVert^2 dx &= \int_{\lvert x \rvert \le 1} \lVert f(x+h) - f(x) \rVert^2 dx +  \int_{\lvert x \rvert > 1} \lVert f(x+h) - f(x) \rVert^2 dx \\
    &\le 2C \cdot \lvert h \rvert^{2 \beta} + C^{1/2}  \cdot \lvert h \rvert^{\beta} \cdot 2 \lVert f \rVert_1  = O(\lvert h \rvert^\beta),
\end{align*}
where $C $ is a positive constant.
Hence, the integral $\int_{\R} \int_{\R} \frac{\lVert f(x+h) - f(x) \rVert^2}{\lvert h \rvert^{1+2s}} dx dh$ is convergent if $\beta > 2s$.
\smallskip
We thus get the following analogue of Corollary \ref{cormesuresstablesinvholdercashilbert} for semigroups.

\smallskip
\bpr \label{propchampsholderiansemigroupesd'op}
Let $H$ be an infinite-dimensional separable complex Hilbert space and let $(T_t)_{t \geq 0}$ be a $\mathcal{C}_0$-semigroup of operators on $H$. Let $\alpha \in (1,2]$. Let $A$ be the generator of this semigroup. Suppose that there exists a finite or countably infinite sequence $(F_j)_{j \in I}$ of maps in $L^1(\R;H) \cap L^2(\R;H)$ such that $A F_j(t) = i t F_j(t)$ for every $j \in I$ and $t \in \R$, and that the subspace $\textrm{span}\{ F_j(t) : j \in I, \, t \in \R \setminus B \}$ is dense in $H$ for every Borel subset $B$ of $\R$ with Lebesgue measure zero.
Suppose also that for each $j \in I$, the map $F_j$ is $\beta_j$-Hölderian for some $\beta_j > 2/\alpha- 1$.
Then the semigroup $(T_t)_{t \geq 0}$ admits a strongly mixing isotropic $\alpha$-stable invariant measure with full support.
\epr

\smallskip
\bpf
The proof follows from the previous remarks, by taking $s_j >0$ such that $\beta_j/2 >  s_j > 1/\alpha -1/2$. Since $ \beta_j >2 s_j$, we obtain that $\int_{\R} \lVert \hat{F_j}(x) \rVert^\alpha dx < \infty$.
\epf

Let us note that Proposition \ref{propchampsholderiansemigroupesd'op} requires to have $\alpha \in (1,2]$, which is less good than Corollary \ref{cormesuresstablesinvholdercashilbert}, but this was expected since we integrate on $\R$. To finish, we have the following result which follows from Plancherel's theorem.

\smallskip
\bpr \label{propchampfourierC^1semigroupes}
Let $H$ be an infinite-dimensional separable complex Hilbert space and let $(T_t)_{t \geq 0}$ be a $\mathcal{C}_0$-semigroup of operators on $H$. Let $\alpha \in (2/3,2]$. Let $A$ be the generator of this semigroup. Suppose that there exists a finite or countably infinite sequence $(F_j)_{j \in I}$ of maps in $L^1(\R;H) \cap L^2(\R;H)$ such that $A F_j(t) = i t F_j(t)$ for every $j \in I$ and $t \in \R$, and that the subspace $\textrm{span}\{ F_j(t) : j \in I, \, t \in \R \setminus B \}$ is dense in $H$ for every Borel subset $B$ of $\R$ with Lebesgue measure zero.
Suppose also that for each $j \in I$, the map $F_j$ is $\mathcal{C}^1$ and that $F_j' $ belongs to $L^1(\R;H) \cap L^2(\R;H)$.
Then the semigroup $(T_t)_{t \geq 0}$ admits a strongly mixing isotropic $\alpha$-stable invariant measure with full support.
\epr

\smallskip

\bpf
Let $f : \R \to H$ be a $\mathcal{C}^1$ function in $L^1(\R; H)\cap L^2(\R;H)$ such that $f'$ also belongs to $L^1(\R; H)\cap L^2(\R;H)$. Then $\widehat{f'}(x) = i x \hat{f}(x)$ on $\R$. Thus, by Plancherel's theorem, we have for $s \in (0,1)$ that
\begin{align*}
    \int_{\R} (1 + \lvert \xi \rvert^2)^s \lVert \hat{f}(\xi) \rVert^2 d\xi &\le \int_{\R} (1 + \lvert \xi \rvert^2) \lVert \hat{f}(\xi) \rVert^2 d\xi \\
    &\lesssim \lVert f \rVert_2^2 + \int_{\R} \lVert \widehat{f'}(\xi) \rVert^2 d\xi \\
    &\lesssim \lVert f \rVert_2^2 + \lVert f' \rVert_2^2.
\end{align*}
Hence $f$ belongs to $\mathcal{H}^{s}(\R;H)$. Since $1/\alpha -1/2 \in [0,1)$, we can just take $s \in (1/\alpha -1/2, 1)$ to obtain that $\int_{\R} \lVert \hat{f}(x) \rVert^\alpha dx < \infty$, concluding the proof of Proposition \ref{propchampfourierC^1semigroupes}.
\epf

\subsection{Examples}  We give two examples of strongly continuous semigroups of operators admitting a strongly mixing isotropic stable invariant measure with full support. We first start with a simple example with weighted shifts. It complements \cite[Thm. 3.11]{Gill1}.

\smallskip

\bpr \label{examplesemigroupshiftspoidsstable}
Let $E = \ell_p(\Z_+)$ with $1 \leq p < \infty$, let $w = (w_n)_{n \geq 1}$ be a sequence of nonzero complex numbers and let $A = B_w$. Set $\beta_0 := 1$ and $\beta_n = w_1 \dotsm w_n$ for $n \geq 1$. Let $T_t = e^{tA}$ for every $t \geq 0$ and suppose that $\rho := \displaystyle \limsup_{n \to \infty} (1/\lvert \beta_n \rvert)^{1/n} < \infty$. Then, for every $\eta \in (0, 1/\rho)$, for every $0 < \alpha <p$ if $p <2$ and $0 < \alpha \le 2$ if $p \ge 2$, the semigroup $(T_t)_{t \geq 0}$ has a strongly mixing isotropic $\alpha$-stable invariant measure with full support.
\epr

\smallskip

\bpf
Let $v_\lambda := (\lambda^n/\beta_n)_{n \geq 0}$ for $\lvert \lambda \rvert <1/\rho := R$ and let $\eta \in (0, 1/\rho)$. Then $A v_\lambda = \lambda v_\lambda$ and so if $F(t) = v_{it}$ for $\lvert t \rvert<\eta$, we obtain that $A F(t) = it F(t)$. Note that the map $t \mapsto F(t) $ if well-defined on $[-\eta, \eta]$ since the convergence radius of $\displaystyle \sum_{n \geq 0} \frac{\lambda^n}{\beta_n}$ is $R$. Since the map $t \mapsto F(t)$ is $\mathcal{C}^\infty$ on $[-\eta, \eta]$, we have that $\lVert \hat{F}(t) \rVert = O(1/\lvert t \rvert^k)$ as $\lvert t \rvert \to \infty$ for every integer $k \ge 1$, by \cite[Proof of Lemma 2.4.3]{HytonenvanNeerverVeraarWeis16}. Finally, it was proved in \cite[Thm. 4.4]{AGNEESSENS_2026} that the subspace $\textrm{span}\{ F(t) : t \in [-\eta, \eta] \setminus B \}$ is dense in $E$ for every Borel subset $B$ of $[-\eta, \eta]$ with Lebesgue measure zero. This concludes the proof of Proposition \ref{examplesemigroupshiftspoidsstable}.
\epf

\smallskip
We now consider translation semigroups. Let $I$ be one of the two intervals $[0,\infty)$ or $\R$, let $p \in [1,\infty)$ and let $\rho : I \to (0,\infty)$ be an admissible weight, i.e., there exist constants $M \ge 1$ and $w \in \R$ such that $\rho(s) \le M e^{wt} \rho(t+s) $ for every $s \in I$ and $t \ge 0$. We consider the Banach space $L_{\rho}^p (I) := \{ f : I \to \C : \int_I \lvert f(x) \rvert^p \rho(x) dx < \infty \}$ equipped with the norm $\lVert f \rVert:= (\int_I \lvert f(x) \rvert^p \rho(x) dx)^{1/p}$. Note that this space has type $\min(2,p)$. The translation semigroup $(T_t)_{t \ge 0}$ is defined on this space by $T_t f = f(\cdot+t)$.

\smallskip

\bpr \label{stableinvariantmeasuretranslationsemigps}
Let $p \in [1,\infty)$, let $I \in \{ [0, \infty), \R \}$, and let $\rho : I \to (0,\infty)$ be an integrable admissible weight.
Suppose that there exists an integer $k \ge 1$ such that $\int_I \lvert x \rvert^{kp} \rho(x) dx < \infty$.
Let $\alpha \in (1/k,2]$ if $p \ge 2$ and $\alpha \in (1/k,p)$ if $p \in [1,2)$. Then, the translation semigroup $(T_t)_{t \ge 0}$ on $L^p_{\rho}(I)$ admits a strongly mixing isotropic $\alpha$-stable invariant measure with full support. 
\epr

\smallskip

\bpf
Let us consider a function $f : \R \to \C$ in $L^1(\R; \C) $ such that $f \ne 0$ almost everywhere and such that $\hat{f}$ belongs to the Schwartz space, and let $e_t : x \in I \mapsto e^{itx}$ for $t \in \R$. Denoting by $F : t \in \R \mapsto f(t) e_t$, we have that $F(t) \in L^p_{\rho}(I)$ for every $t \in \R$, and $T_t F(s) = e^{ist} F(s)$ for every $t \ge 0$ and $s \in \R$. Moreover, let us notice that $F \in L^1(\R;L^p_{\rho}(I))$. Indeed, we have that
\begin{align*}
    \int_{\R} \lVert F(s) \rVert ds &= \int_{\R} \left (\int_{I} \lvert f(s) \rvert^p \rho(x) dx \right)^{1/p} ds \\
    &= \int_{\R} \lvert f(s) \rvert \left( \int_{I} \rho(x) dx\right)^{1/p} ds <\infty.
\end{align*}
Moreover, let us notice that 
\begin{align*}
    \int_{\R} \lVert \hat{F}(x) \rVert^\alpha dx = \int_{\R} \left( \int_I \lvert \hat{f}(x-y) \rvert^p \rho(y) dy \right)^{\alpha/p} dx.
\end{align*}
Let $C>0$ be a constant such that $\lvert \hat{f}(t) \rvert^p \le C (1 + \lvert t \rvert)^{-kp}$ for every $t \in \R$. Then
\begin{align*}
    \int_{I} \lvert \hat{f}(x-y) \rvert^p \rho(y) dy &\le C \int_I (1 + \lvert x-y \rvert)^{-kp} \rho(y) dy \\ 
    &\le (1 + \lvert x \rvert)^{-kp } \int_I (1 + \lvert y \rvert)^{kp} \rho(y) dy.
\end{align*}
Hence
\begin{align*}
    \int_{\R} \left( \int_{I} \lvert \hat{f}(x-y) \rvert^p \rho(y) dy \right)^{\alpha/p} dx \le \left( \int_{\R} (1+ \lvert x \rvert)^{-k \alpha} dx \right) \cdot \left( \int_I (1 + \lvert y \rvert)^{kp} \rho(y) dy \right)^{\alpha/p} < \infty 
\end{align*}
since $k \alpha >1$, since $\rho$ is integrable and since $\int_I \lvert y \rvert^{kp} \rho(y) dy < \infty$.
Finally, since $f \ne 0$ almost everywhere, the subspace $\textrm{span}\{ F(t) : t \in \R \setminus B \}$ is dense in $L^2_{\rho}(I)$ for every Borel subset $B$ of $\R$ with Lebesgue measure zero, by \cite[Proof of Thm. 4.2]{AGNEESSENS_2026}. This concludes the proof of Proposition \ref{stableinvariantmeasuretranslationsemigps} by Theorem \ref{thminvmeasurestablesemigroups}. 
\epf

\smallskip
Proposition \ref{stableinvariantmeasuretranslationsemigps} is the same result as \cite[Thm. 3.10]{Gill1}, but with a different proof. In the case where $p=2$ and $\alpha = 2$, the conditions in Proposition \ref{stableinvariantmeasuretranslationsemigps} require that $\int_I \lvert x \rvert^{2k} \rho(x) dx < \infty$ for some integer $k \ge 1$. In \cite[p. 1139]{AGNEESSENS_2026}, it is required to have $\int_I \lvert x \rvert^{4} \rho(x) dx < \infty$. In fact, for $p = 2$ and $\alpha = 2$, we do not need more assumption than the integrability of the admissible weight $\rho$. 

\smallskip

\bth \label{Gaussianinvmeasuresemigrouptranslations}
Let $p = 2$, let $I \in \{ [0, \infty), \R \}$, and let $\rho : I \to (0,\infty)$ be an integrable admissible weight. Then, the translation semigroup $(T_t)_{t \ge 0}$ on $L^2_\rho(I)$ admits a strongly mixing isotropic Gaussian invariant measure with full support.
\eth

\smallskip

\bpf
Let $f \in L^1(\R;\C) \cap L^2(\R; \C)$ such that $f \ne 0$ almost everywhere. Since $\rho$ and $f$ are integrable, the function $F : t \in \R \mapsto f(t) e_t$ is integrable on $\R$. Moreover, we also have that
\begin{align*}
    \int_{\R} \lVert F(s) \rVert^2 ds = \int_{\R} \left ( \int_I \rho(x) dx \right) \lvert f(s) \rvert^2 ds < \infty.
\end{align*}
The conclusion follows from Proposition \ref{propexistencemesureinvgaussiennesemigroupe}.
\epf

\smallskip
Finally, for the case $p = 2$, Proposition \ref{stableinvariantmeasuretranslationsemigps} shows that if $\int_I \lvert x \rvert^{2} \rho(x) dx < \infty$, then we obtain a strongly mixing $\alpha$-stable invariant measure for every $\alpha \in (1,2]$. We improve a bit the lower bound by $\alpha \in (2/3,2]$ under this integrability condition on the weight.  

\smallskip

\bpr
Let $p = 2$, let $I \in \{ [0, \infty), \R \}$, and let $\rho : I \to (0,\infty)$ be an integrable admissible weight such that $\int_I \lvert x \rvert^2 \rho(x) dx < \infty$. Then, for every $\alpha \in (2/3,2]$, the translation semigroup $(T_t)_{t \ge 0}$ on $L^2_\rho(I)$ admits a strongly mixing isotropic $\alpha$-stable invariant measure with full support.
\epr

\smallskip

\bpf
Let $f(x) = e^{-x^2}$ for every $x \in \R$. Let again $F : t \in \R \mapsto f(t) e_t$, which belongs to $L^1(\R;L^2_\rho(I)) \cap L^2(\R; L^2_\rho(I))$. Then the map $F$ is $\mathcal{C}^1$ with $F'(t)(x) = f'(t) e^{itx} + f(t) i x e^{itx}$. We can easily check that 
\begin{align*}
    \int_{\R} \lVert F'(t) \rVert^2 dt \le \int_{\R} \lvert f'(t) \rvert^2 \left (\int_I \rho(x) dx \right) dt \, + \int_{\R} \lvert f(t) \rvert^2 \left (\int_I x^2 \rho(x) dx \right) dt < \infty. 
\end{align*}
And using that $(a+b)^{1/2} \le a^{1/2} + b^{1/2}$ for every $a,b \ge 0$, we can also show that
\begin{align*}
    \int_{\R} \lVert F(t) \rVert dt \le \int_{\R} \lvert f'(t) \rvert \left ( \int_I \rho(x) dx \right)^{1/2} dt + \int_{\R} \lvert f(t) \rvert \left(\int_I x^2 \rho(x) dx \right)^{1/2} dt < \infty. 
\end{align*}
The conclusion follows from Proposition \ref{propchampfourierC^1semigroupes}.
\epf

\section{Additional results and open questions}

We end this article with some comments and open questions in relation to our previous results.
\smallskip

First, let us notice that our construction of the random integral with respect to a complex infinitely divisible random measure on $(J,\Sigma,m)$ requires to have that the function $f : J \to E$ is in $L^p(J;E)$, when the Banach space is of type $p \in [1,2]$. In fact, we have the following result. 

\smallskip
\bpr \label{resultintegralealeatoiremesureID}
Let $E$ be a separable complex Banach space of type $p \in [1,2]$ and let $M$ be a complex infinitely random measure on a $\sigma$-finite measure space $(J, \Sigma,m)$ satisfying $\E(M(A)) = 0$ and $\E(\lvert M(A) \rvert^p) \le C' m(A)$ for every measurable set $A$ of finite measure, for some constant $C'>0$. Then $\int_J f dM$ exists if $\int_J \min(1, \lVert f(t) \rVert^p) m(dt) < \infty$.
\epr

\smallskip
Of course, the integrability condition in Proposition \ref{resultintegralealeatoiremesureID} is the same as $\int_J \lVert f(t) \rVert^p m(dt)$ in the case where $J$ is a subset of $\Z$ equipped with the counting measure. But in the case where $J$ is an interval of $\R$ with the Lebesgue measure, it may be weaker than $\int_J \lVert f(t) \rVert^p m(dt) < \infty$, depending on the function $f$ and on the interval $J$.
\smallskip

\bpf[Proof of Proposition \ref{resultintegralealeatoiremesureID} ]
We rely on \cite[pp. 190--191]{Woyczynski19}. Let us write that $\int_J f dM = \int_J f \indic_{\{ \lVert f \rVert \le 1 \}} dM + \int_J f \indic_{\{ \lVert f \rVert >1 \}} dM$, where the first random integral is well defined under the assumption  $\int_J \min(1, \lVert f(t) \rVert^p) m(dt) < \infty$, by the classical construction of the random integral. For the second random integral, we will show that the series $\int_J f \indic_{\{ \lVert f \rVert >1 \}} dM = \displaystyle \sum_{k \ge 1} \int_{k < \lVert f \rVert \le k+1} f dM$ is almost surely convergent. Note that the random integrals in the sums are well defined too. Let $q_k := \frac{\log(2)}{\log(k+1)}$ for every $k \ge 1$. Since $E$ is of type $q$ for any $q \in (0,p]$, we have, by (\ref{inegalitéesperanceBanachtypep}), that
\begin{align*}
    \E \left( \left \lVert \int_{k < \lVert f \rVert \le k+1} f dM \right \rVert^{q_k} \right) &\le C \int_{k < \lVert f \rVert \le k+1} \lVert f(t) \rVert^{q_k} m(dt) \\
    &\leq 2C m(\{ k < \lVert f \rVert \le k+1 \})
\end{align*}
for every $k \ge 1$. By Markov's inequality we obtain that
\begin{align*}
    \sum_{k \ge 1} \P \left( \left \lVert \int_{k < \lVert f \rVert \le k+1} f dM \right\rVert > k^{-2} \right) &\le 2C \sum_{k \ge 1} \frac{m(\{ k < \lVert f \rVert \le k+1 \})}{k^{-2 q_k}} \\
    &\le C_1 m(\{ \lVert f \rVert > 1 \}) \\
    &< \infty,
\end{align*}
where $C_1$ is a positive constant given by the convergence of $k^{-2q_k}$ to $e^{-2 \log(2)}$. An application of the Borel-Cantelli theorem concludes the proof of Proposition \ref{resultintegralealeatoiremesureID}.
\epf

Let us also notice that, in the proof of Proposition \ref{mesurestableinvshiftspoids}, we had to redefine the random integral for the case $p>2$ in order to make use of the summability assumption on the weight. The problem was that the Banach space $E = \ell_p(\Z_+)$ is of type $2$, but we did not necessarily have that $\displaystyle \sum_{n \ge 0} \lVert  \hat{F}(n) \rVert^2 < \infty$. More generally, our construction of the random integral $\int_J f dM$ with respect to complex infinitely divisible random measures requires to have the $p$-integrability of the function $f$ in the case where the Banach space is of type $p$, which can be very restrictive. We believe that the construction made by Rosinski in \cite{Rosinski87} may be more practical for our problem. More precisely, motivated by the proof of Proposition \ref{mesurestableinvshiftspoids} and by Rosinski \cite{Rosinski87}, we say that a measurable function $f : J \to E$ is integrable with respect to a complex infinitely divisible random measure $M$ if there exists a sequence of simple functions $(f_n)_{n \ge 1}$ converging almost everywhere to $f$ such that for every measurable set $A \in \Sigma$, the sequence $(\int_A f_n dM)_{n \ge 1}$ converges in probability. In this case, we set $\int_A f dM := \displaystyle \lim_{n \to \infty} \int_A f_n dM$ for every $A \in \Sigma$, where the convergence holds in probability. The random integral $\int_J f dM$ is well defined by \cite[Lemma 3.1.2]{Rosinski87}, and the characteristic function of the random integral is given as in the introduction of this paper, by \cite[Lemma 3.1.4]{Rosinski87}. Moreover, the map $A \mapsto \int_A f dM$ is an $E$-valued independently scattered random measure on $(J, \Sigma)$. We believe that with this construction of the random integral, we can recover the result of Bayart and Grivaux \cite[Thm 5.38]{BM09}, and we leave it as an \textbf{open problem}.

\medskip
Finally, we would like to come back to Kalisch-type operators. We showed in Proposition \ref{Kalisch-typestableLP} that for every $\alpha \in (1,2]$, the operator $T$ admits a strongly mixing isotropic $\alpha$-stable invariant measure with full support. The conclusion followed from the inequality $1/p > 1/\alpha - 1/p'$, together with the fact that the Fourier coefficients of the unimodular eigenvector field are $\alpha$-summable if $\alpha >1$. We claim that, for $\alpha \in (0,1]$, the operator $T$ admits an isotropic $\alpha$-stable invariant measure with full support. Indeed, let $(y_n)_{n \geq 1}$ be a dense sequence in $[0,2\pi)$ and let $x_n := 2^{-n} \indic_{(y_n, 2\pi)}$ for each $n \ge 1$. Then we have that $T x_n = \lambda_n x_n$ with $\lambda_n := e^{i y_n}$ for every $n \ge 1$, and $\displaystyle \sum_{n \ge 1} \lVert x_n \rVert^\alpha< \infty$. Moreover, the subspace $\textrm{span}\{ x_n : n \ge 1 \}$ is dense in $E = L^p[0,2\pi]$ for $1<p< \infty $. Indeed, if $g \in L^{p'}[0,2\pi]$ satisfies $\int_{y_n}^{2\pi} g(t) dt = 0$ for each $n \ge 1$, then $G = 0$, where $G(y) := \int_y^{2\pi} g(t) dt$, since $g$ is integrable and $G$ is continuous on $ [0,2\pi)$. This implies that $g = 0$ almost everywhere, proving the density result. The existence of an isotropic $\alpha$-stable invariant measure with full support for the operator $T$ thus follows from \cite[Proposition 4.5]{Gill1}. However, this gives a non-ergodic stable invariant measure. Moreover, let us notice that the proof given in Proposition \ref{Kalisch-typestableC0} does not work for $\alpha \in (0,1]$. Hence, the question of whether the operator $T$ in Propositions \ref{Kalisch-typestableLP} and \ref{Kalisch-typestableC0} admits an ergodic $\alpha$-stable invariant measure for $\alpha \in (0,1]$ is left as an \textbf{open question}.

\bibliography{bibdossier}

@book {GEP,
	AUTHOR = {Grosse-Erdmann, Karl-G. and Peris Manguillot, Alfredo},
	TITLE = {Linear chaos},
	SERIES = {Universitext},
	PUBLISHER = {Springer, London},
	YEAR = {2011},
	PAGES = {xii+386},
	ISBN = {978-1-4471-2169-5},
	MRCLASS = {47-02 (37D45 47A16)},
	MRNUMBER = {2919812},
	MRREVIEWER = {E.\ A.\ Gallardo-Guti\'errez},
	DOI = {10.1007/978-1-4471-2170-1},
	URL = {https://doi.org/10.1007/978-1-4471-2170-1},
}

@book {BM09,
    AUTHOR = {Bayart, Fr\'ed\'eric and Matheron, \'Etienne},
     TITLE = {Dynamics of linear operators},
    SERIES = {Cambridge Tracts in Mathematics},
    VOLUME = {179},
 PUBLISHER = {Cambridge University Press, Cambridge},
      YEAR = {2009},
     PAGES = {xiv+337},
      ISBN = {978-0-521-51496-5},
   MRCLASS = {47-02 (11M06 37B05 47A16 47A35)},
  MRNUMBER = {2533318},
MRREVIEWER = {E.\ A.\ Gallardo-Guti\'errez},
       DOI = {10.1017/CBO9780511581113},
       URL = {https://doi.org/10.1017/CBO9780511581113},
}

@article {BayartGrivaux06,
    AUTHOR = {Bayart, Fr\'ed\'eric and Grivaux, Sophie},
     TITLE = {Frequently hypercyclic operators},
   JOURNAL = {Trans. Amer. Math. Soc.},
  FJOURNAL = {Transactions of the American Mathematical Society},
    VOLUME = {358},
      YEAR = {2006},
    NUMBER = {11},
     PAGES = {5083--5117},
      ISSN = {0002-9947,1088-6850},
   MRCLASS = {47A16 (37A05 37B05 47A35)},
  MRNUMBER = {2231886},
MRREVIEWER = {Vivien\ G.\ Miller},
       DOI = {10.1090/S0002-9947-06-04019-0},
       URL = {https://doi.org/10.1090/S0002-9947-06-04019-0},
}

@article {BayartGrivaux07,
    AUTHOR = {Bayart, Fr\'ed\'eric and Grivaux, Sophie},
     TITLE = {Invariant {G}aussian measures for operators on {B}anach spaces
              and linear dynamics},
   JOURNAL = {Proc. Lond. Math. Soc. (3)},
  FJOURNAL = {Proceedings of the London Mathematical Society. Third Series},
    VOLUME = {94},
      YEAR = {2007},
    NUMBER = {1},
     PAGES = {181--210},
      ISSN = {0024-6115,1460-244X},
   MRCLASS = {47A16 (28C20 28D05 37A30 47A35)},
  MRNUMBER = {2294994},
MRREVIEWER = {Douglas\ P.\ Dokken},
       DOI = {10.1112/plms/pdl013},
       URL = {https://doi.org/10.1112/plms/pdl013},
}

@article {ChakirMourchid18,
    AUTHOR = {Chakir, M. and EL Mourchid, S.},
     TITLE = {Strong mixing {G}aussian measures for chaotic semigroups},
   JOURNAL = {J. Math. Anal. Appl.},
  FJOURNAL = {Journal of Mathematical Analysis and Applications},
    VOLUME = {459},
      YEAR = {2018},
    NUMBER = {2},
     PAGES = {778--788},
      ISSN = {0022-247X,1096-0813},
   MRCLASS = {47D07 (28C20 28D05 60G15)},
  MRNUMBER = {3732554},
MRREVIEWER = {Marta\ Tyran-Kami\'nska},
       DOI = {10.1016/j.jmaa.2017.11.003},
       URL = {https://doi.org/10.1016/j.jmaa.2017.11.003},
}

@article {MauPrivault,
    AUTHOR = {Mau, Camille and Privault, Nicolas},
     TITLE = {Mixing of linear operators under infinitely divisible measures
              on {B}anach spaces},
   JOURNAL = {J. Math. Anal. Appl.},
  FJOURNAL = {Journal of Mathematical Analysis and Applications},
    VOLUME = {535},
      YEAR = {2024},
    NUMBER = {1},
     PAGES = {Paper No. 128160, 26},
      ISSN = {0022-247X,1096-0813},
   MRCLASS = {60B11 (28C20 47A35 60E07)},
  MRNUMBER = {4699523},
       DOI = {10.1016/j.jmaa.2024.128160},
       URL = {https://doi.org/10.1016/j.jmaa.2024.128160},
}

@article{AGNEESSENS_2026, title={Frequently hypercyclic random vectors for $\mathcal{C}_{0}$-semigroups}, DOI={10.1017/etds.2026.10274}, journal={Ergodic Theory and Dynamical Systems}, author={Agneessens, KEVIN}, year={2026}, pages={1–30}}

@article {AGNEESSENS_2023,
    AUTHOR = {Agneessens, Kevin},
     TITLE = {Frequently hypercyclic random vectors},
   JOURNAL = {Proc. Amer. Math. Soc.},
  FJOURNAL = {Proceedings of the American Mathematical Society},
    VOLUME = {151},
      YEAR = {2023},
    NUMBER = {3},
     PAGES = {1103--1117},
      ISSN = {0002-9939,1088-6826},
   MRCLASS = {47A16 (47A35 47B37 47B38)},
  MRNUMBER = {4531641},
       DOI = {10.1090/proc/16153},
       URL = {https://doi.org/10.1090/proc/16153},
}

@article {MouzeMunnier21,
    AUTHOR = {Mouze, Augustin and Munnier, Vincent},
     TITLE = {Frequent hypercyclicity of random holomorphic functions for
              {T}aylor shifts and optimal growth},
   JOURNAL = {J. Anal. Math.},
  FJOURNAL = {Journal d'Analyse Math\'ematique},
    VOLUME = {143},
      YEAR = {2021},
    NUMBER = {2},
     PAGES = {615--637},
      ISSN = {0021-7670,1565-8538},
   MRCLASS = {47A16 (30H99 60F15 60G50)},
  MRNUMBER = {4299171},
       DOI = {10.1007/s11854-021-0164-5},
       URL = {https://doi.org/10.1007/s11854-021-0164-5},
}

@article {BayartMatheron16,
    AUTHOR = {Bayart, Fr\'ed\'eric and Matheron, \'Etienne},
     TITLE = {Mixing operators and small subsets of the circle},
   JOURNAL = {J. Reine Angew. Math.},
  FJOURNAL = {Journal f\"ur die Reine und Angewandte Mathematik. [Crelle's
              Journal]},
    VOLUME = {715},
      YEAR = {2016},
     PAGES = {75--123},
      ISSN = {0075-4102,1435-5345},
   MRCLASS = {47A35 (46B09 47C15)},
  MRNUMBER = {3507920},
MRREVIEWER = {Makoto\ Mori},
       DOI = {10.1515/crelle-2014-0002},
       URL = {https://doi.org/10.1515/crelle-2014-0002},
}

@book {SamorodnitskyTaqqu94,
    AUTHOR = {Samorodnitsky, Gennady and Taqqu, Murad S.},
     TITLE = {Stable non-{G}aussian random processes},
    SERIES = {Stochastic Modeling},
      NOTE = {Stochastic models with infinite variance},
 PUBLISHER = {Chapman \& Hall, New York},
      YEAR = {1994},
     PAGES = {xxii+632},
      ISBN = {0-412-05171-0},
   MRCLASS = {60E07 (60G18 60G99 60H99 62J05 62M10 62M15)},
  MRNUMBER = {1280932},
MRREVIEWER = {Thomas\ Mikosch},
}

@article {Grivaux11,
    AUTHOR = {Grivaux, Sophie},
     TITLE = {A new class of frequently hypercyclic operators},
   JOURNAL = {Indiana Univ. Math. J.},
  FJOURNAL = {Indiana University Mathematics Journal},
    VOLUME = {60},
      YEAR = {2011},
    NUMBER = {4},
     PAGES = {1177--1201},
      ISSN = {0022-2518,1943-5258},
   MRCLASS = {47A16 (37A05 46B09 47A35)},
  MRNUMBER = {2975340},
MRREVIEWER = {Fr\'ed\'eric\ Bayart},
       DOI = {10.1512/iumj.2011.60.4350},
       URL = {https://doi.org/10.1512/iumj.2011.60.4350},
}

@article {MurilloPeris15,
    AUTHOR = {Murillo-Arcila, M. and Peris, A.},
     TITLE = {Strong mixing measures for {$C_0$}-semigroups},
   JOURNAL = {Rev. R. Acad. Cienc. Exactas F\'is. Nat. Ser. A Mat. RACSAM},
  FJOURNAL = {Revista de la Real Academia de Ciencias Exactas, F\'isicas y
              Naturales. Serie A. Matematicas. RACSAM},
    VOLUME = {109},
      YEAR = {2015},
    NUMBER = {1},
     PAGES = {101--115},
      ISSN = {1578-7303,1579-1505},
   MRCLASS = {47D06 (35L20)},
  MRNUMBER = {3315705},
MRREVIEWER = {Rodica\ Luca},
       DOI = {10.1007/s13398-014-0169-3},
       URL = {https://doi.org/10.1007/s13398-014-0169-3},
}

@article {Flytzanis95,
    AUTHOR = {Flytzanis, E.},
     TITLE = {Unimodular eigenvalues and linear chaos in {H}ilbert spaces},
   JOURNAL = {Geom. Funct. Anal.},
  FJOURNAL = {Geometric and Functional Analysis},
    VOLUME = {5},
      YEAR = {1995},
    NUMBER = {1},
     PAGES = {1--13},
      ISSN = {1016-443X,1420-8970},
   MRCLASS = {28D05 (47A35 60B11)},
  MRNUMBER = {1312018},
MRREVIEWER = {Wojciech\ Bartoszek},
       DOI = {10.1007/BF01928214},
       URL = {https://doi.org/10.1007/BF01928214},
}

@article {MurilloPeris13,
    AUTHOR = {Murillo-Arcila, M. and Peris, A.},
     TITLE = {Strong mixing measures for linear operators and frequent
              hypercyclicity},
   JOURNAL = {J. Math. Anal. Appl.},
  FJOURNAL = {Journal of Mathematical Analysis and Applications},
    VOLUME = {398},
      YEAR = {2013},
    NUMBER = {2},
     PAGES = {462--465},
      ISSN = {0022-247X,1096-0813},
   MRCLASS = {47A16},
  MRNUMBER = {2990071},
       DOI = {10.1016/j.jmaa.2012.08.050},
       URL = {https://doi.org/10.1016/j.jmaa.2012.08.050},
}

@book {Woyczynski19,
    AUTHOR = {Woyczy\'nski, Wojbor A.},
     TITLE = {Geometry and martingales in {B}anach spaces},
 PUBLISHER = {CRC Press, Boca Raton, FL},
      YEAR = {2019},
     PAGES = {xiii+315},
      ISBN = {978-1-138-61637-0},
   MRCLASS = {46-02 (46B20 46L53 60B11 60G46)},
  MRNUMBER = {3890006},
MRREVIEWER = {Olav\ Nygaard},
}

@book {HytonenvanNeerverVeraarWeis16,
    AUTHOR = {Hyt\"onen, Tuomas and van Neerven, Jan and Veraar, Mark and
              Weis, Lutz},
     TITLE = {Analysis in {B}anach spaces. {V}ol. {I}. {M}artingales and
              {L}ittlewood-{P}aley theory},
    SERIES = {Ergebnisse der Mathematik und ihrer Grenzgebiete. 3. Folge. A
              Series of Modern Surveys in Mathematics [Results in
              Mathematics and Related Areas. 3rd Series. A Series of Modern
              Surveys in Mathematics]},
    VOLUME = {63},
 PUBLISHER = {Springer, Cham},
      YEAR = {2016},
     PAGES = {xvi+614},
      ISBN = {978-3-319-48519-5; 978-3-319-48520-1},
   MRCLASS = {46-02 (42B35 46E30)},
  MRNUMBER = {3617205},
MRREVIEWER = {Adam\ Os\polhk ekowski},
}

@book {LiQueffelec18,
    AUTHOR = {Li, Daniel and Queff\'elec, Herv\'e},
     TITLE = {Introduction to {B}anach spaces: analysis and probability.
              {V}ol. 1},
    SERIES = {Cambridge Studies in Advanced Mathematics},
    VOLUME = {166},
      NOTE = {Translated from the French by Dani\`ele Gibbons and Greg
              Gibbons,
              For the French original see [MR2124356]},
 PUBLISHER = {Cambridge University Press, Cambridge},
      YEAR = {2018},
     PAGES = {xxx+431},
      ISBN = {978-1-107-16051-4; 978-1-107-16263-1},
   MRCLASS = {46B03 (46Bxx 60B12 60G15)},
  MRNUMBER = {3729312},
MRREVIEWER = {Ramon\ van Handel},
}

@book {Linde86,
    AUTHOR = {Linde, Werner},
     TITLE = {Probability in {B}anach spaces---stable and infinitely
              divisible distributions},
    SERIES = {A Wiley-Interscience Publication},
   EDITION = {Second},
 PUBLISHER = {John Wiley \& Sons, Ltd., Chichester},
      YEAR = {1986},
     PAGES = {195},
      ISBN = {0-471-90893-2},
   MRCLASS = {60B11 (28B05 46G12 47B10 60E07)},
  MRNUMBER = {874529},
}

@book {Schwartz81,
    AUTHOR = {Schwartz, Laurent},
     TITLE = {Geometry and probability in {B}anach spaces},
    SERIES = {Lecture Notes in Mathematics},
    VOLUME = {852},
      NOTE = {Based on notes taken by Paul R. Chernoff},
 PUBLISHER = {Springer-Verlag, Berlin-New York},
      YEAR = {1981},
     PAGES = {x+101},
      ISBN = {3-540-10691-X},
   MRCLASS = {46-02 (46Bxx 60B11)},
  MRNUMBER = {612976},
MRREVIEWER = {Niels\ J\o rgen\ Nielsen},
}

@book {EngelNagel06,
    AUTHOR = {Engel, Klaus-Jochen and Nagel, Rainer},
     TITLE = {A short course on operator semigroups},
    SERIES = {Universitext},
 PUBLISHER = {Springer, New York},
      YEAR = {2006},
     PAGES = {x+247},
      ISBN = {978-0387-31341-2; 0-387-31341-9},
   MRCLASS = {47-01 (47D03 47D06)},
  MRNUMBER = {2229872},
MRREVIEWER = {Jacek\ Banasiak},
}

@article {Rosinski87,
    AUTHOR = {Rosi\'nski, Jan},
     TITLE = {Bilinear random integrals},
   JOURNAL = {Dissertationes Math. (Rozprawy Mat.)},
  FJOURNAL = {Polska Akademia Nauk. Instytut Matematyczny. Dissertationes
              Mathematicae. Rozprawy Matematyczne},
    VOLUME = {259},
      YEAR = {1987},
     PAGES = {71},
      ISSN = {0012-3862},
   MRCLASS = {60H05 (46G10 60B11 60E07)},
  MRNUMBER = {888463},
MRREVIEWER = {Terry\ R.\ McConnell},
}

@article{Gill1,
  author = {Valentin Gillet},
  title  = {Stable invariant measures in linear dynamics},
  status = {preprint},
  eprint = {arxiv:2606.30463},
  date   = {2026},
}

@article {Konig91,
    AUTHOR = {K\"onig, Hermann},
     TITLE = {On the {F}ourier-coefficients of vector-valued functions},
   JOURNAL = {Math. Nachr.},
  FJOURNAL = {Mathematische Nachrichten},
    VOLUME = {152},
      YEAR = {1991},
     PAGES = {215--227},
      ISSN = {0025-584X,1522-2616},
   MRCLASS = {46E40 (42B05)},
  MRNUMBER = {1121235},
MRREVIEWER = {Raymond\ H.\ Cox},
       DOI = {10.1002/mana.19911520118},
       URL = {https://doi.org/10.1002/mana.19911520118},
}

@article {BayartGrivauxMatheronMenet25,
    AUTHOR = {Bayart, Fr\'ed\'eric and Grivaux, Sophie and Matheron,
              \'Etienne and Menet, Quentin},
     TITLE = {Hereditarily frequently hypercyclic operators and disjoint
              frequent hypercyclicity},
   JOURNAL = {Ergodic Theory Dynam. Systems},
  FJOURNAL = {Ergodic Theory and Dynamical Systems},
    VOLUME = {45},
      YEAR = {2025},
    NUMBER = {10},
     PAGES = {3021--3072},
      ISSN = {0143-3857,1469-4417},
   MRCLASS = {47A16 (37A30 37B20)},
  MRNUMBER = {4955406},
MRREVIEWER = {Khadija\ Elkhalloufy},
       DOI = {10.1017/etds.2025.13},
       URL = {https://doi.org/10.1017/etds.2025.13},
}

@book {Hormander90,
    AUTHOR = {H\"ormander, Lars},
     TITLE = {The analysis of linear partial differential operators. {I}},
    SERIES = {Grundlehren der mathematischen Wissenschaften [Fundamental
              Principles of Mathematical Sciences]},
    VOLUME = {256},
   EDITION = {Second},
      NOTE = {Distribution theory and Fourier analysis},
 PUBLISHER = {Springer-Verlag, Berlin},
      YEAR = {1990},
     PAGES = {xii+440},
      ISBN = {3-540-52345-6},
   MRCLASS = {35-02 (42B10 46Fxx)},
  MRNUMBER = {1065993},
       DOI = {10.1007/978-3-642-61497-2},
       URL = {https://doi.org/10.1007/978-3-642-61497-2},
}

@article {BenyiOh13,
    AUTHOR = {B\'enyi, \'Arp\'ad and Oh, Tadahiro},
     TITLE = {The {S}obolev inequality on the torus revisited},
   JOURNAL = {Publ. Math. Debrecen},
  FJOURNAL = {Publicationes Mathematicae Debrecen},
    VOLUME = {83},
      YEAR = {2013},
    NUMBER = {3},
     PAGES = {359--374},
      ISSN = {0033-3883,2064-2849},
   MRCLASS = {42B35 (42B10 42B25)},
  MRNUMBER = {3119672},
MRREVIEWER = {Anna\ K.\ Savvopoulou},
       DOI = {10.5486/PMD.2013.5529},
       URL = {https://doi.org/10.5486/PMD.2013.5529},
}

@article {ImekrazRobertThomann16,
    AUTHOR = {Imekraz, Rafik and Robert, Didier and Thomann, Laurent},
     TITLE = {On random {H}ermite series},
   JOURNAL = {Trans. Amer. Math. Soc.},
  FJOURNAL = {Transactions of the American Mathematical Society},
    VOLUME = {368},
      YEAR = {2016},
    NUMBER = {4},
     PAGES = {2763--2792},
      ISSN = {0002-9947,1088-6850},
   MRCLASS = {60G50 (42A32)},
  MRNUMBER = {3449257},
MRREVIEWER = {Mikhail\ A.\ Lifshits},
       DOI = {10.1090/tran/6607},
       URL = {https://doi.org/10.1090/tran/6607},
}

@book {DiestelJarchowTonge95,
    AUTHOR = {Diestel, Joe and Jarchow, Hans and Tonge, Andrew},
     TITLE = {Absolutely summing operators},
    SERIES = {Cambridge Studies in Advanced Mathematics},
    VOLUME = {43},
 PUBLISHER = {Cambridge University Press, Cambridge},
      YEAR = {1995},
     PAGES = {xvi+474},
      ISBN = {0-521-43168-9},
   MRCLASS = {46-02 (46Bxx 46M05 47-02 47B10 47D50)},
  MRNUMBER = {1342297},
MRREVIEWER = {Andreas\ Defant},
       DOI = {10.1017/CBO9780511526138},
       URL = {https://doi.org/10.1017/CBO9780511526138},
}

@article {MarcusPisier84,
    AUTHOR = {Marcus, M. B. and Pisier, G.},
     TITLE = {Characterizations of almost surely continuous {$p$}-stable
              random {F}ourier series and strongly stationary processes},
   JOURNAL = {Acta Math.},
  FJOURNAL = {Acta Mathematica},
    VOLUME = {152},
      YEAR = {1984},
    NUMBER = {3-4},
     PAGES = {245--301},
      ISSN = {0001-5962,1871-2509},
   MRCLASS = {60G17 (60E07)},
  MRNUMBER = {741056},
MRREVIEWER = {R.\ M.\ Dudley},
       DOI = {10.1007/BF02392199},
       URL = {https://doi.org/10.1007/BF02392199},
}

@article {Fouriersurvey98,
    AUTHOR = {Garsia-Kuerva, Kh. and Kazaryan, K. S. and Kolyada, V. I. and
              Torrea, Kh.\ L.},
     TITLE = {The {H}ausdorff-{Y}oung inequality with vector-valued
              coefficients and applications},
   JOURNAL = {Uspekhi Mat. Nauk},
  FJOURNAL = {Uspekhi Matematicheskikh Nauk},
    VOLUME = {53},
      YEAR = {1998},
    NUMBER = {3(321)},
     PAGES = {3--84},
      ISSN = {0042-1316,2305-2872},
   MRCLASS = {42B10 (42A99 42C99 46B20 46E40)},
  MRNUMBER = {1657592},
MRREVIEWER = {Plamen\ Dzhakov},
       DOI = {10.1070/rm1998v053n03ABEH000018},
       URL = {https://doi.org/10.1070/rm1998v053n03ABEH000018},
}

@book {VianaOliveira16,
    AUTHOR = {Viana, Marcelo and Oliveira, Krerley},
     TITLE = {Foundations of ergodic theory},
    SERIES = {Cambridge Studies in Advanced Mathematics},
    VOLUME = {151},
 PUBLISHER = {Cambridge University Press, Cambridge},
      YEAR = {2016},
     PAGES = {xvi+530},
      ISBN = {978-1-107-12696-1},
   MRCLASS = {37-02 (28D05 37Axx 37D35 54H20)},
  MRNUMBER = {3558990},
MRREVIEWER = {Douglas\ P.\ Dokken},
       DOI = {10.1017/CBO9781316422601},
       URL = {https://doi.org/10.1017/CBO9781316422601},
}

\end{document}